\documentclass[12pt]{amsart}
\usepackage{amsmath, amssymb, amsthm, amsfonts, mathrsfs}
\usepackage{cite}
\usepackage{amscd}
\usepackage{url}

\usepackage{graphicx}
\usepackage{caption}
\usepackage{subcaption}
\usepackage{fullpage}
\usepackage{lscape, pdflscape}

\theoremstyle{definition}

\begin{document}
\title{Inferring Fitness in Finite Populations with Moran-like dynamics}

\author{Marc Harper}
\address{University of California Los Angeles}
\email{marcharper@ucla.edu} 
\date{\today}

\begin{abstract}
Biological fitness is not an observable quantity and must be inferred from population dynamics. Bayesian inference applied to the Moran process and variants yields a robust inference method that can infer fitness in populations evolving via a Moran dynamic and generalizations. Information about fitness is derived solely from birth-events in birth-death and death-birth processes in which selection acts proportionally to fitness, which allows the method to be applied to populations on a network where the network itself may be changing in time. Populations may also be allowed to change size while still allowing estimates for fitness to be inferred.
\end{abstract}

\maketitle

\section{Introduction}

The theory of natural selection dominates biological explanation and thinking. Despite the elegance of Darwin's approach to descent with modification and later refinements and abstractions, accurate models and measurements of real populations of replicating organisms can be very difficult to obtain. Analogies of fitness and physical potentials notwithstanding, fitness is more like a hidden variable than an observable quantity. There is no instrument that measures fitness in a manner analogous to a voltmeter. Determining a hidden quantity can be achieved with statistical inference. In this paper we discuss such a method in finite populations with geometric structure, and address critical theoretical questions such as: How much information does a single birth or death event yield about the fitness of a replicating organism?

Statistical inference has been successfully applied in many areas of biology and bioinformatics, e.g. in phylogenetics \cite{felsenstein2004inferring}. The approach in this paper uses models that differ from some other recent approaches to inferring fitness, such as \cite{shaw2010inferring} and \cite{otwinowski2013genotype}, yet the goal is much the same: practically quantify evolutionary parameters such as fitness that are very useful intuitively. In particular, we will be concerned with traditional models from population biology, like the Moran process\cite{moran1962statistical} \cite{moran1962statistical} and recent extensions that include graphical distribution \cite{lieberman2005evolutionary}. This is to bring together the knowledge gained from evolutionary game theory, population dynamics, and statistical inference rather than a statement on which models are more appropriate. Throughout this investigation, it will become clear that separating birth and death events yields greater information gain per population state transition, so we will consider some new variants of the Moran model. We study the effectiveness of Bayesian inference applied to these models by simulating large numbers of population trajectories.

Let us first consider a basic version of the general problem of inferring fitness. Consider a population of two types of replicating entities $A$ and $B$, with $a$ individuals of type $A$ and $b$ individuals of type $B$, with a total finite population size of $N = a + b$. Type can refer to genotype or phenotype as is appropriate. Suppose that replicators of type $A$ have fitness $f_A = 1$ and replicators of type $B$ have $f_B = r$, where $r = f_B / f_A$ is the relative fitness between the two types. Given that the population is in the state $(a, b)$ where both $a > 0$ and $b > 0$ (so that both types of replicators are represented in the population), let an individual be selected for replication proportionally due to fitness. More precisely, an individual of type $A$ is selected with probability $a / (a + r b)$ and an individual of type $B$ is selected with probability $ r b / (a + rb)$. The state of the population $(a, b)$ is not generally informative about the relative fitness $r$ because it is not possible to know without the history of the population what the relative differences imply about $r$. Concretely, if $a > b$, if could be that type $A$ is dominating the population and $b$ is declining, or that type $B$ has a relatively new phenotype that is currently in the process of fixating and has yet to take over the majority of the population. If it were known that the state was stable under some suitable definition of stability then $r$ could be estimated by $a / b$ (assuming $a \approx rb$, which occurs when $a \approx rN / (r + 1)$). This, however, assumes many things that are likely to be unrealistic in a variety of scenarios -- for example (1) that the state of the population can be known precisely (in the laboratory perhaps, but less likely for natural populations), and (2) that the population has converged to a stable state (which implies that the population has been under observation for some time, or was in a stable configuration upon observation). Ideally we would prefer a method that allows one to come upon a population, observe it for a short time, capturing possibly only a random sample of selection events and population states, and make an accurate estimate of the relative fitness $r$. This is achievable.

One way to approach the problem of inferring $r$ would be to use the fixation probabilities of type $B$ for the population in configuration $(a, b)$, which is given by $\phi = \frac{1 - r^{-b}}{1 - r^{-N}}$ \cite{nowak2006evolutionary}. Then by preparing many populations and recording the proportion of populations that fixation with all members of type $A$ or type $B$, the ratio of such events could be used to solve numerically for $r$, by inferring $\phi$ as a binomial parameter. The approach has several flaws: (1) it is likely costly and very time consuming, if even practically possible; (2) it is not applicable to populations on graphical structures for which fixation formulas are not known and may be practically impossible to calculate; and (3) it ignores all the information revealed by the population trajectory -- that is, all the information obtainable from replication events -- using only the information obtainable from the final state.

What does an observation of a replication of either type $A$ or $B$ tell us about the relative fitness $r$? Let us suppose first that $a \approx b$. A replication of an individual type $B$ then suggests that $r > 1$. How much information such an event gives about about $r$ depends on the relative sizes of the quantities $a$ and $rb$. If $a$ is much larger than $b$, a replication of an individual of type $B$ is evidence that $r$ is larger than one (since based on $a$ and $b$ alone, a replication of a type $A$ replicator should have been more likely), whereas if $a \approx rb$ then the selection landscape was approximately uniform, and a replication of type $A$ was just as likely. In this case, the replication event yields little information about the relative fitness. Of course, we cannot put too much stock into a single event (in case of a fluke). Repeated observations strengthen the information yield.

Now let us make the claim that replication events yield information more precise. Assume that we can keep the population in the state $(a,b)$, where $N = a + b$, and observe arbitrarily many replication events. Let $\alpha$ and $\beta$ count the number of observed replications of type $A$ and $B$ respectively. We can then estimate the value of $a / (a + rb)$ with a maximum likelihood estimate $\alpha / (\alpha + \beta)$, and using the values of $a$ and $b$ determine $r$. This is simply a version of Bayes' original method of inference for determining the value of a binomial parameter $\theta$ using Beta distributions with parameters $\alpha$ and $\beta$, with $\theta = a / (a + rb) \in [0,1]$, and distribution function:
\[ f(\theta; \alpha, \beta) = \frac{\Gamma(\alpha + \beta)}{\Gamma(\alpha) \Gamma(\beta)} \theta^{\alpha - 1} (1 - \theta)^{\beta - 1}, \]
where $\Gamma$ is the generalized factorial function with $\Gamma(n) = (n - 1)!$ for integers $n$. The estimate $\hat{\theta} = \alpha / (\alpha + \beta)$ is the mean of the Beta distribution. To quantify how much information is gained about $r$ from each individual replication event, fix $\alpha$ and $\beta$ and observe one birth event. This yields $\alpha'$ and $\beta'$ where either $\alpha' = \alpha + 1$ or $\beta' = \beta + 1$. Then the information gain associated to the birth event is given by the Kullback-Leibler information divergence  $D_{KL}(\text{Beta}(\alpha', \beta') || \text{Beta}(\alpha, \beta))$ \cite{kullback1951information}. To be clear, this is the information gained experimentally by the observer about $r$, not the information gained by the population about its environment. Information about the relative fitness is inferred from the population trajectory (in this case, the sequence of states of the population had we not artificially fixed the state).

Unlike in this example, a population will vary among a sequence of population states as birth and death events occur, possibly changing in total size and changing its spatial distribution. At each state, the probabilities of each type replicating will differ since they depend on the population state $(a, b)$. We no longer have a single $\theta$ (which depends on $r$) to infer and it will be more direct to work with probability distributions for $r$ on $[0, \infty)$ or $[0, R]$ for some maximum value $R$. It is also no longer the case that a sample of identically distributed observations is possible -- observations of birth and death events will depend on the state and change the state as a population evolves (altering the fitness landscape). Nevertheless, the conclusions of some typical results for similar inference problems hold in simulations despite the lack of these common assumptions.

The Moran process is a birth-death process that models natural selection in finite populations \cite{moran1958random} \cite{moran1962statistical}. Study of the Moran process has been extensive, including fixation probabilities for various landscapes and starting states \cite{fudenberg2004stochastic} \cite{taylor2004evolutionary} \cite{antal2006fixation} \cite{nowak2006evolutionary}, evolutionary stability \cite{schaffer1988evolutionarily} \cite{fogel1998instability} \cite{ficici2000effects}, the evolution of cooperation \cite{nowak2004emergence}, and in the context of multi-level selection \cite{traulsen2005stochastic}. Unlike deterministic models of selection such as the replicator equation \cite{hofbauer2003evolutionary}, the Moran process is a stochastic process, specifically a Markov chain. A Markov process does not have a fixed trajectory for a given starting state, so when we refer to a trajectory we mean one particular sequence of states of a population modeled by the process. Markov chains are typically analyzed in terms of transition probabilities and quantities derivable from the transitions, such as absorption probabilities and mean convergence times. Such analysis for the Moran process can be found in many sources, such as \cite{nowak2006evolutionary}. For our purposes, we will need actual population sequences so as to infer values of $r$ from population trajectories. We will consider the Moran process and variations that separate birth and death events, including generalizations of the Moran process to populations distributed on graphs, which are well known in evolutionary graph theory \cite{lieberman2005evolutionary} \cite{pacheco2006active}. In particular, we will consider birth-death processes on graphs in which a member of the population is selection for reproduction from the entire population, replacing an outbound neighbor, and death-birth processes in which a member is randomly selected to be replaced by an inbound neighbor (which is chosen proportionally to fitness).

\section{Preliminaries}

\subsection{Markov Processes for Birth and Death}
The Moran process models selection in a well-mixed finite population of replicating entities. As before, consider a population of $N$ replicators, $a$ of type A (A-individuals) and $b$ of type B (B-individuals) where $N = a + b$ is fixed as before. Individuals of type A and B have fitness $f_A$ and $f_B$ respectively, and may depend on the population parameters $a$ and $b$. Although we could dispense with one of the parameters $a$ or $b$ since $N$ is fixed, we will need both for a modification that does not maintain a fixed value for $N$. We will only consider populations of integral size $N$ at least 3.

The population is updated by selecting an individual at random proportionally to fitness and selecting an individual at random uniformly to be replaced. For the population to change state, individuals of different types must be chosen. The transition probabilities between states are given by \cite{taylor2004evolutionary} \cite{nowak2004emergence}
\begin{align*}
T_{(a, b) \to (a+1, b-1)} &= \frac{a f_A}{a f_A + b f_B} \frac{b}{a+b} \\
T_{(a, b) \to (a-1, b+1)} &= \frac{b f_B}{a f_A + b f_B} \frac{a}{a+b}, \\
T_{(a, b) \to (a, b)} &= \frac{a^2 f_A + b^2 f_B}{a f_A + b f_B} \frac{1}{a+b}
\end{align*}
with $T_{(a, b) \to (a, b)} = 1 - T_{(a, b) \to (a+1, b-1)} - T_{(a, b) \to (a-1, b+1)}$. The fitness landscape is given by
\begin{align*}
f_A(i) &= \frac{a(i-1) + b(N-i)}{N-1} \\
f_B(i) &= \frac{ci + d(N-i-1)}{N-1}
\end{align*}
for a game matrix defined by
\[ \left( \begin{matrix}
 a & b\\
 c & d
\end{matrix} \right). \]
The Moran process is given by $a=1=c, b=r=d$. The process has two absorbing states $(0, N)$ and $(N, 0)$, corresponding to the fixation of one of the two types. We will first consider the Moran process for populations with $F_A = 1$ and $F_B = r$ where $r \in [0, \infty)$ is the relative fitness of type B versus type A.
Given $r$, much can be said about the Moran process in this case, including fixation probabilities starting from any population state. We consider a different problem. Given a trajectory of the Markov process as a series of states can we accurately determine the value of $r$? We could directly measure the rates of reproduction and compare those values, which we will refer to as \emph{counting}. From the Moran process we can only detect birth of either type if the population changes state; otherwise we could only guess based on the population distribution which type both had a birth and a death. The combined birth-death transition probabilities wash out the information obtained in the case that the population state stays the same (e.g. individuals of the same type replicates and then dies). For this reason we also consider a modification of the process that breaks the transitions into separate birth and death events. The results of the process are the same as the Moran process but we will know in each step which type replicated and which type died. In this case, the probabilities for fitness proportionate reproduction are given by:
\begin{align*}
T_{(a, b) \to (a+1, b)} &= \frac{a f_A}{a f_A + b f_B} \\
T_{(a, b) \to (a, b+1)} &= \frac{b f_B}{a f_A + b f_B} \\
\end{align*}
We can also view this process as conditional on the size of the population. If $a + b = N$, we select a replicator to reproduce. If the population is of size $N+1$, then we randomly select a replicator to be removed. Inference benefits from this modification because the transitions which leave the state unchanged now yield information. Since the death event is random, it carries no information about $r$, only about the population size, which is assumed known in this case.  We could alternatively consider a death-birth process. Results in this case are similar but slightly more susceptible to stochastic noise for very small population states. Because the death event occurs first, the reproductive event is at a state of smaller population size $N-1$, but otherwise basically the same as the birth-death case. For populations on graphs, birth-death and death-birth processes can be quite different.

\subsection{Bayesian Inference}
To use Bayesian inference to infer the value of $r$ from a sequence of states we first choose a prior probability distribution for the possible values of $r$. In theory $r \in [0, \infty)$, but for computational purposes we will restrict to an interval $r \in [0, R]$, where $R$ is chosen to be sufficiently large. Choosing a prior probability distribution is a subjective process, though in this case there are a few reasonable choices. We could use a uniform distribution on $[0, R]$ in attempt to chose an uninformative distribution. The uniform prior, however, puts too much weight on the interval $[1, R]$ which should have approximately the same weight as on $[0,1]$ since $r$ is a ratio. By the neutral theory of evolution \cite{kimura1985neutral}, most mutations will have little effect on the relative fitness $r$, so priors which are somehow centered at $r=1$ are reasonable choices in many scenarios. It is also reasonable to assume that $p(0) = 0 = p(R)$ since both types $A$ and $B$ exist in the population (and presumably can reproduce). A gamma distribution with well-chosen parameters fits both assumptions and can have mean, median, or mode at $r=1$. We will not dwell on the problem of prior choice other than to say that a different prior could further improve the accuracy for some combinations of parameters and states. In simulations, gamma distributions based on the assumptions of the neutral theory perform well. Later in the text we will see that there is another choice for the prior distribution that is computationally convenient. The gamma distribution takes the form
\[ f(r; k, \theta) = \frac{r^{k-1}e^{-x / \theta}}{\theta^k \Gamma(k)}, \]
where $\theta$ is unrelated to the previous usage. Note that there are many distributions that could take the place of the gamma distribution and there is no particular reason for it over other distributions that also have the desired properties, though it is noteworthy that the tail of the distribution drops off exponentially. Simulations reported in this paper use parameters $k=2$ and $\theta=2$.

Bayesian inference takes a distribution on the possible values of the parameter and an observation to produce an updated distribution that takes into account the information gained by the observation. Starting with the prior distribution, inference proceeds until the data is exhausted or the parameter is known to a sufficient accuracy. To update the distribution using the sequence of states we use Bayes' Theorem, where $E$ indicates a transition from $(a, b) \to (a', b')$:
\[ Pr(r | E) = \frac{Pr(E|r) Pr(r)}{P(E)},\]
and we have that $Pr(E|r) = T_{(a, b) \to (a', b')}(r)$ from the definition of the Markov process. $P(E)$ can be calculated by integrating over $r$, but it is only necessary to normalize at the final step unless we wish to have intermediate estimates. To form an estimate over the entire sequence of states, pair the states into transitions $E_1, \ldots , E_n$ (each state is in two adjacent pairs), and form the posterior distribution $\text{Posterior}(r| E_1, \ldots , E_n) \propto \text{Prior}(r) \prod_{i}{Pr(E_i|r)}$. Finally, to extract an estimate for $r$, we normalize the posterior and compute the mean value $\int_{0}^{R}{r \text{Posterior}(r) \, dr}$ or mode numerically.

\subsection{Conjugate distributions}
Just as in the introductory example, the observations of birth events (or population state changes) can be used to determine probability distributions for the purposes of inference. In the case of a binomial distribution, each observation of a success or failure multiplies the inferred distribution by $\theta$ or $1 - \theta$ (increasing the value of $\alpha$ or $\beta$). In the case of a birth-death process as described above, we multiply by the corresponding transition factor. Now $\alpha$ and $\beta$ are vector parameters indexed by $a=1, \ldots, N-1$ (with $b = N-a$) and the distribution takes the form
\[ P(r) \propto \prod_{a=1}^{N-1}{\left( \frac{a}{a + r b} \right)^{\alpha_a} \left( \frac{r b}{a + r b} \right)^{\beta_a}}. \]
Rewriting, and introducing a normalization constant $N_{\alpha, \beta}$ yields the form
\[ P(r) = N_{\alpha, \beta} \prod_{a=1}^{N-1}{\frac{r^{\beta_a}}{(a + r (N-a))^{\alpha_a + \beta_a}}} \]
The normalization constant can be computed symbolically using partial fractions but the formulas are, even for simple values of $\alpha$ and $\beta$, very complex for even relatively small values of $N$; similarly so for the integration to determine the mean of the distribution. Moreover, analytic formulas for the maximum likelihood estimator and the Fisher information are not easily obtainable. For this reason we proceed numerically, noting that using this conjugate distribution still prevents integration at every step in the inference process. We will refer to this distribution as the FPS distribution (fitness proportionate selection distribution).

The FPS distribution provides more options for the prior distribution $P$. If any $\beta_a > 0$ then $P(0) = 0$. Similarly, if any $\alpha_a > 0$ then $P(r) \to 0$ as $r \to \infty$. For the distribution to have a finite integral in the case that the distribution is supported on $[0, \infty)$, we need that $\sum_a \alpha_a > 1$. So, one possible choice of parameters for a prior is $\alpha_a = 1 = \beta_a$ for all $a$. For these parameters, the FPS distribution has mode equal to 1 (the maximum likelihood estimate) and mean approximately $1.16$. The mode is equal to one if $\alpha_a = \beta_{N-a}$ for all $a$, and more generally if $a \beta_a  = \alpha_a (N-a)$ for all $a$. To see this, notice that the equation for derivative of the logarithm of $P$ to be zero can be rewritten as
\[ r = \frac{\sum_{a=1}^{N-1}{ \prod_{i \neq a}{(i + r (N - i))} a \beta_a}}{\sum_{a=1}^{N-1}{ \prod_{i \neq a}{(i + r (N - i))} (N-a) \alpha_a}}.\]
For $r$ close to one, we have the approximation $r \approx \sum_{a=1}^{N-1}{a \beta_a} / \sum_{a=1}^{N-1}{(N-a) \alpha_a}$. Merely counting the replication events for each type would lead to an estimate of the form $r = \sum_{a=1}^{N-1}{\beta_a} / \sum_{a=1}^{N-1}{\alpha_a}$, so it is clear how the position in the population affects the estimate (as well as the true value of $r$). For this reason we will modify the counting procedure in a later section to weight the value of a replication event by the population state, which significantly improves the results. In general, the maximum likelihood method is a given by the real positive root of a polynomial and is unique in $[0, \infty)$, which can be shown with basic techniques from calculus. Given a prior, the estimate for $r$ given by the mode of the posterior is known as the maximum a posteriori probability estimate (equal to the maximum likelihood estimate if the prior is uniform). This estimate can be easily computed numerically. Note that the full distribution can be used to give a credible interval in addition to a point estimate. 

For the Moran process, the analogous distribution includes a third vector parameter $\gamma$ corresponding to transitions in which the population state is unchanged:
\[ P(r) \propto \prod_{a=1}^{N-1}{\left( \frac{a}{a + r b} \frac{b}{a + b} \right)^{\alpha_a} \left( \frac{r b}{a + r b} \frac{a}{a+b} \right)^{\beta_a}  \left( \frac{a^2 + r b^2}{a + r b} \frac{1}{a+b} \right)^{\gamma_a}}. \]
Similarly to the previous case, we can rewrite this as
\[ P(r) \propto \prod_{a=1}^{N-1}{ \frac{r^{\beta_a} \left(a^2 + r (N-a)^2 \right)^{\gamma_a}}{ (a + r (N-a))^{\alpha_a + \beta_a + \gamma_a}}} \]

\begin{figure}
    \centering
    \includegraphics[scale=0.3]{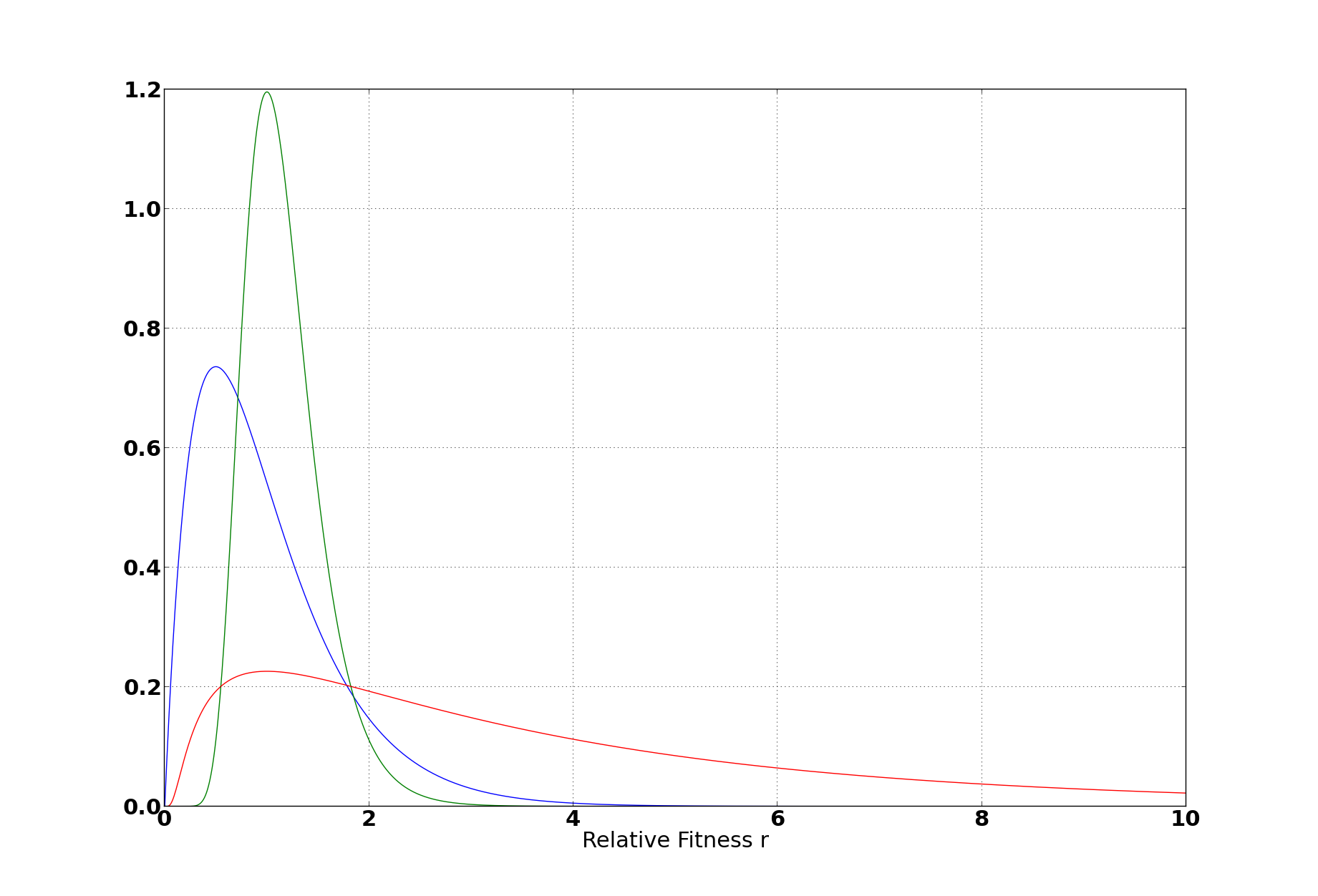}
    \caption{Examples of Priors. Blue: Gamma with $k=2=\theta$. Green: FPS distribution with N=30, $\alpha_a = 1 = \beta_a$ for all $a$. Red: FPS with $\beta_1 = 10 = \alpha_1$ and all others 0.}
    \label{priors}
\end{figure}

\subsection{Conjugate Distributions for Variably-Sized Populations}

We will also consider the case that the population size may change. To do so, we merely need to include terms for different values of $N$, where $\alpha$ and $\beta$ are now triangular matrix parameters $\alpha_{N, a}$ and $\beta_{N, a}$. The FPS distribution is given by
\begin{align*}
P(r) &\propto \prod_{N \in \mathcal{N}}\prod_{a=1}^{N-1}{\left( \frac{a}{a + r b} \right)^{\alpha_{N, a}} \left( \frac{r b}{a + r b} \right)^{\beta_{N, a}}}\\ 
&\propto \prod_{N \in \mathcal{N}} \prod_{a=1}^{N-1}{\frac{r^{\beta_{N,a}}}{(a + r (N-a))^{\alpha_{N,a} + \beta_{N,a}}}}
\end{align*}
where $\mathcal{N}$ is the possible set of total population sizes for the particular process. Call this distribution the variable-population FPS distribution. The previously described FPS distribution is the special case when $\mathcal{N} = \{ N \}$. This distribution would be used, for instance, with a Markov process in which, rather than have a birth and death event in each cycle, has a birth or death event with some probability $p$ that potentially depends on the population state (but not the parameter to be inferred). In this case multiple birth or death events could occur in sequence without the other; still only the birth events will be used to infer the value of $r$. In particular, suppose a population has carrying capacity $N = 2K$, so that the probability of a death event is $p=1$ (so that the population cannot exceed its carrying capacity), and where $p < 1$ for all population states in which $N = a + b < 2K$, perhaps given by a sigmoid function with inflection point at $N=K$.

From these examples it should be clear that any Markov chain on population states with parameter dependent transition probabilities could be treated in a similar manner. This paper will only cover the distributions described so far, but the reader could consider more complex examples where the fitness functions and transition probabilities depend on more variables (such as mutation probabilities or fitness functions deriving from unknown game matrices).

\subsection{Conjugate Distributions for Populations on Graphs}

Suppose now that rather than a well-mixed population, we have populations in which the replicators occupy vertices on a directed graph, such as a directed cycle or a graph with an undirected star topology. For the birth-death process on a graph, a replicator is chosen proportionally to fitness and replaces a randomly selected outgoing neighbor. Since the replicator that reproduces is selected from the entire population, the appropriate probability distribution is the fixed-population FPS distribution. If the graph can change size (i.e. lose or gain vertices), then the variable-population FPS distribution is needed.

For death-birth processes, the variable-population FPS distribution is needed for both static graphs and those in which vertices and edges may be added or removed. In this case, the set $\mathcal{N}$ consists of the sizes of the sets of inbound neighbors for each vertex (since these are the subpopulations that in which replicators are selected for reproduction).

\section{Computations and Simulations}
Trajectories for the Markov process were computed with \emph{mpsim} (\textbf{M}arkov \textbf{p}rocess \textbf{sim}ulator), a flexible and parallelized simulator for any Markov process for which the transition graph can be specified and stored in computer memory, created specifically for the purpose of generating large numbers of trajectories. Source code for \emph{mpsim} is available on Github at \url{https://github.com/marcharper/mpsim}. Additional code containing more simulations and code to process trajectories, compute posterior probability distributions and relative fitness estimates, is also available on Github at \url{https://github.com/marcharper/fitness_inference}. Parameter estimates are ultimately determined by numerical integration.

\begin{figure}
    \centering
    \includegraphics[width=\textwidth]{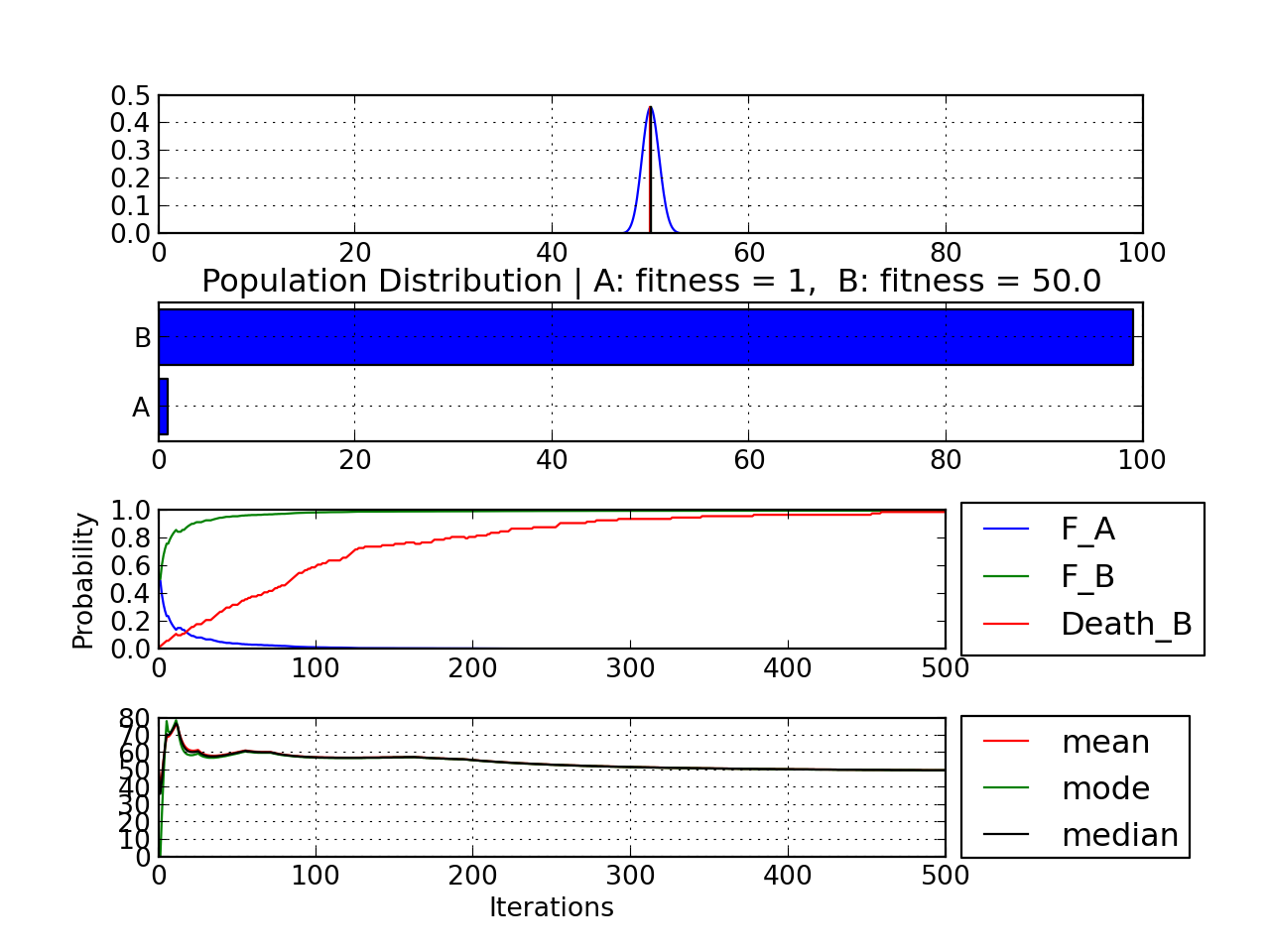}
    \caption{The End Result of Inference. In this diagram, the population began at $a=1, b=99$. The top plot shows the posterior distribution after 500 transitions (iterations). The third plot shows the probabilities of selecting type $A$ and $B$ to reproduce, $F_A$ and $F_B$ respectively, as well as the probability of choosing a type $B$ individual to be replaced ($\text{Death}_B$). The lowest plot shows the mean, median, and mode of the posterior distribution over all 500 iterations.}
    \label{final_distribution}
\end{figure}

\begin{figure}
    \centering
    \includegraphics[scale=0.3]{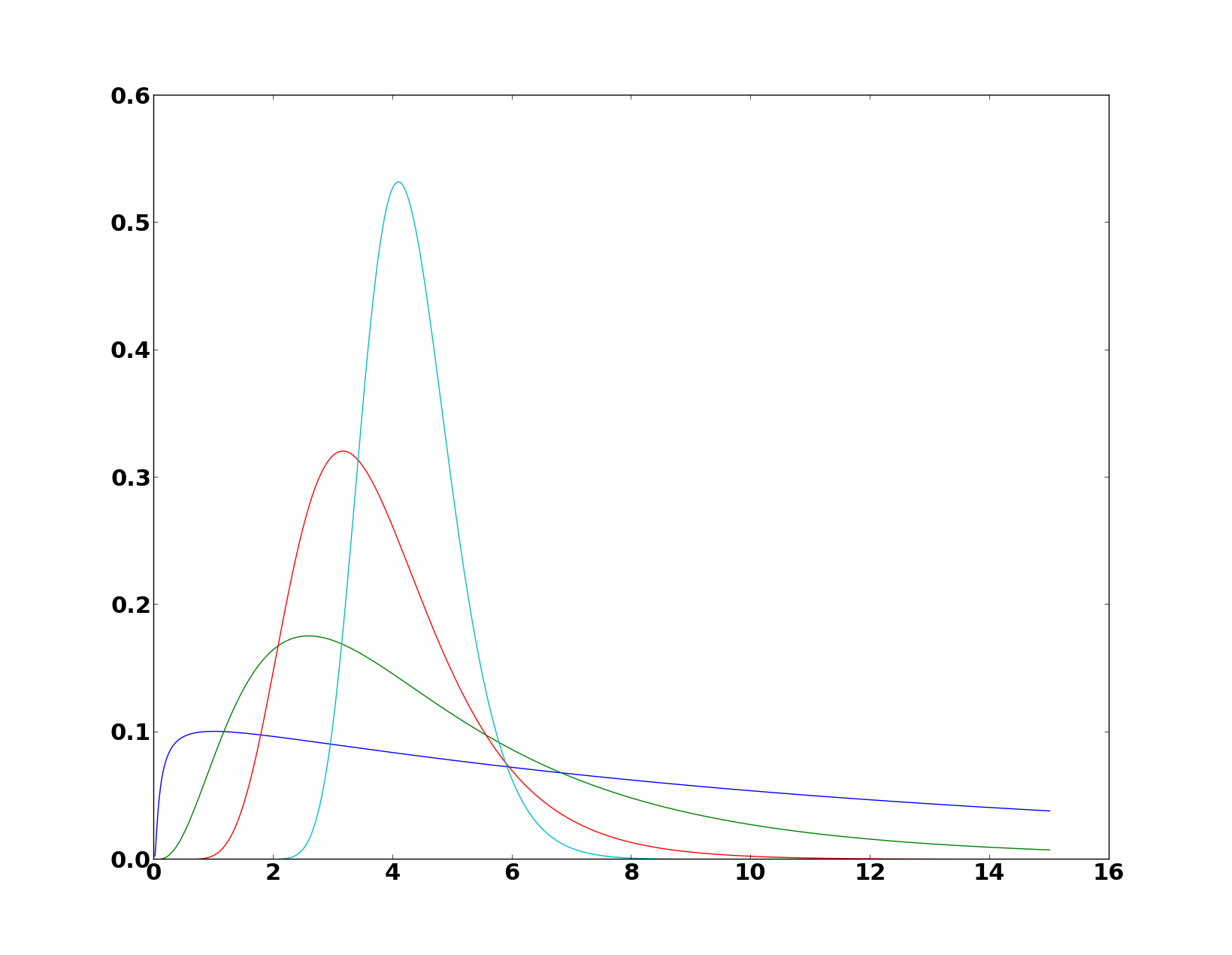}
    \caption{PDFs for various stages of one example of inference. Blue: Prior (FPS with $\alpha_1 = 5 = \beta_1$); Green: After 15 steps; Red: After 60 steps; Cyan: Final posterior (920 steps). Population size $N=60$, true value of $r$ is 4; the inferred value is $r=4.3$.}
    \label{example_inference_run}
\end{figure}

\begin{figure}
    \centering
    \includegraphics[scale=0.3]{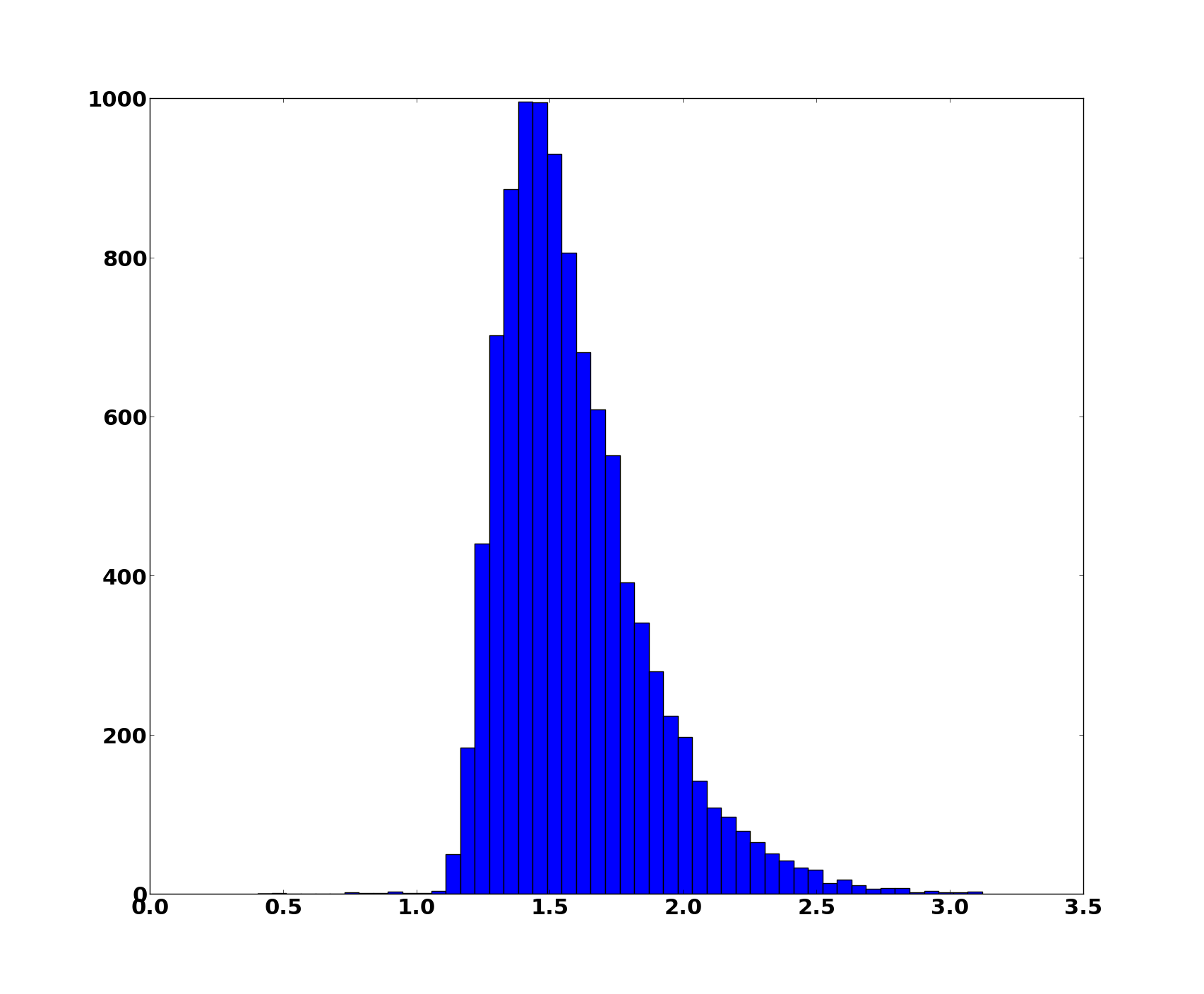}
    \caption{Histogram of inferred final values for 10,000 simulations with $N=40$, $r=1.5$, starting state $(24, 16)$ for the Moran process. The prior distribution is a Gamma with $k=2=\theta$. A small portion of trajectories estimate $r < 1$.}
    \label{moran_histogram}
\end{figure}

\subsection{Simulation Results, Moran Process}
For the purposes of comparison, for each trajectory the total counts of observations of birth events of the two types is recorded to yield estimates of the form
\[ \hat{r} = \frac{\sum_a{\beta_a}}{\sum_a{\alpha_a}},\]
mimicking the maximum likelihood method as discussed in previous sections. These estimates are compared to numerically computed means and modes of posterior distributions.

Estimates for the value of $r$ by counting and by inference via simulation reveal several commonalities. 1000 sample runs for $N$ ranging from 3 to 50, for $r$ ranging from 0.1 to 2 in steps of 0.1, and for each possible starting state were generated and the methods of determining $r$ compared. It is not possible to present results from all possible starting states for all combinations of variables, so we will focus on the most interesting initial starting points: a single mutant, and an initial state with fitness proportionate selection probabilities equal (or as nearly so as possible). The full data set is available on request (and can be regenerated with the previously referred to software). See Figures \ref{moran_sims_balance} and \ref{moran_sims_single_mutant}.

\begin{landscape}
\centering
\begin{figure}[h]
        \begin{subfigure}[b]{0.4\textwidth}
            \centering
            \includegraphics[width=\textwidth]{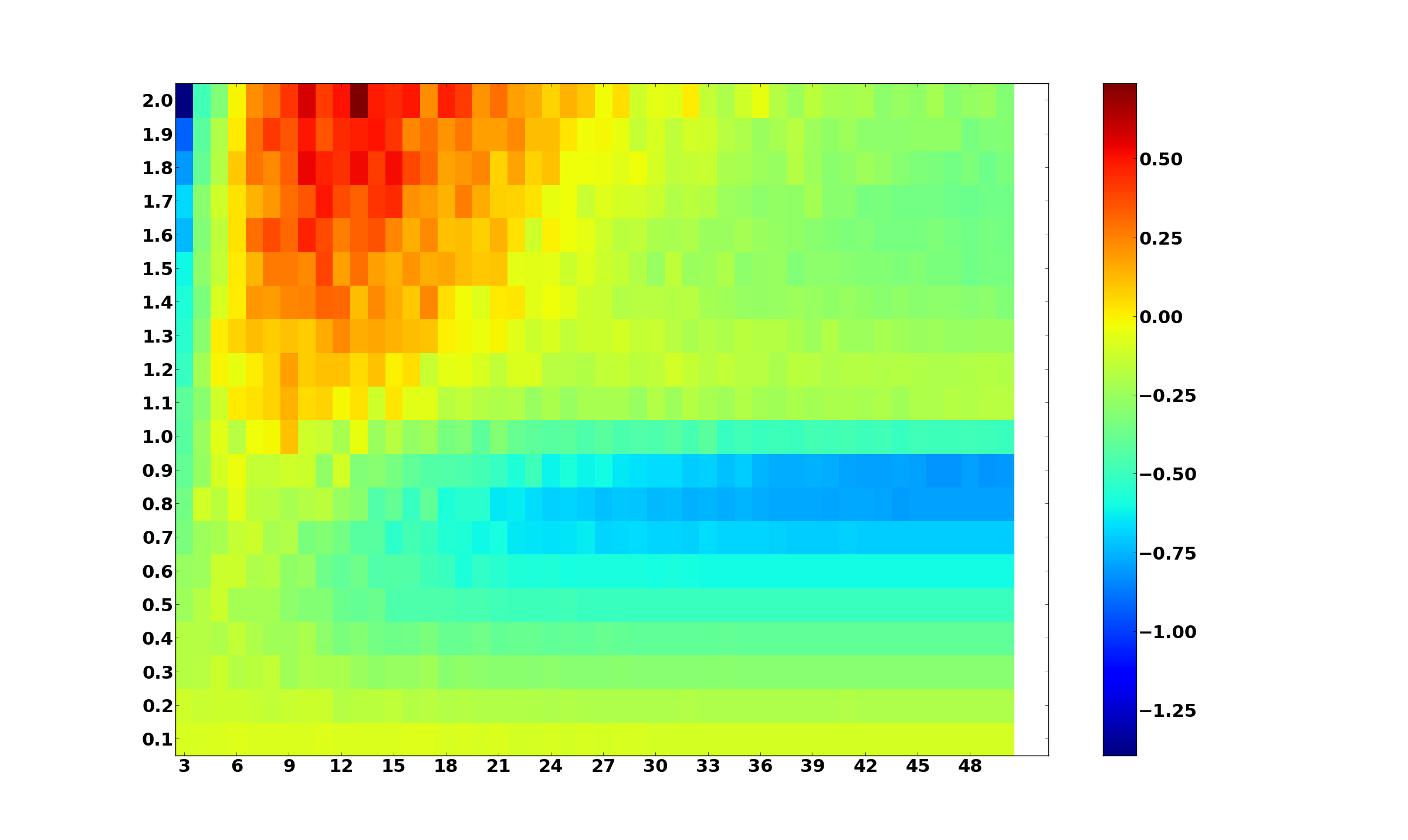}
        \end{subfigure}%
        ~ 
        \begin{subfigure}[b]{0.4\textwidth}
            \centering
            \includegraphics[width=\textwidth]{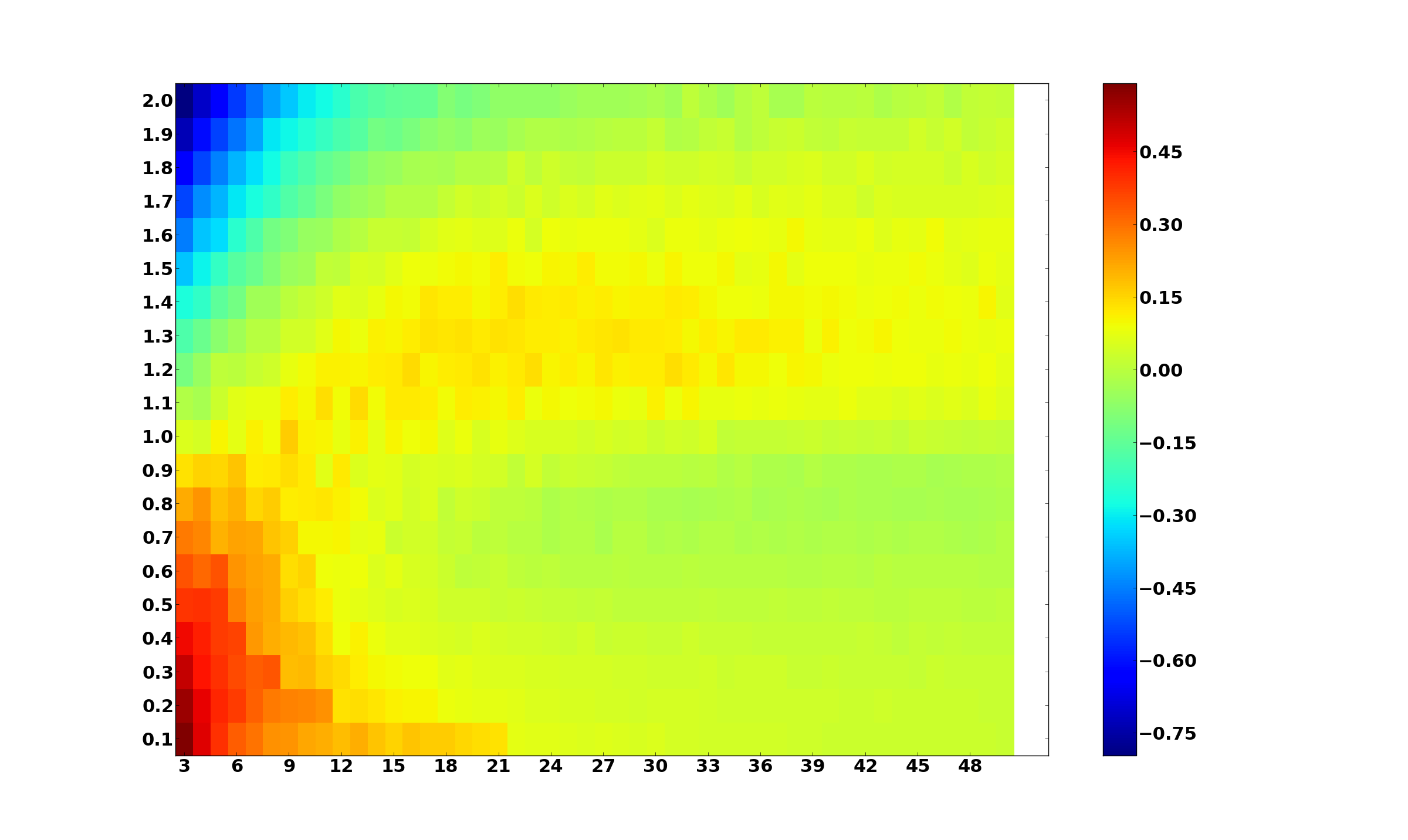}
        \end{subfigure}
        \begin{subfigure}[b]{0.4\textwidth}
            \centering
            \includegraphics[width=\textwidth]{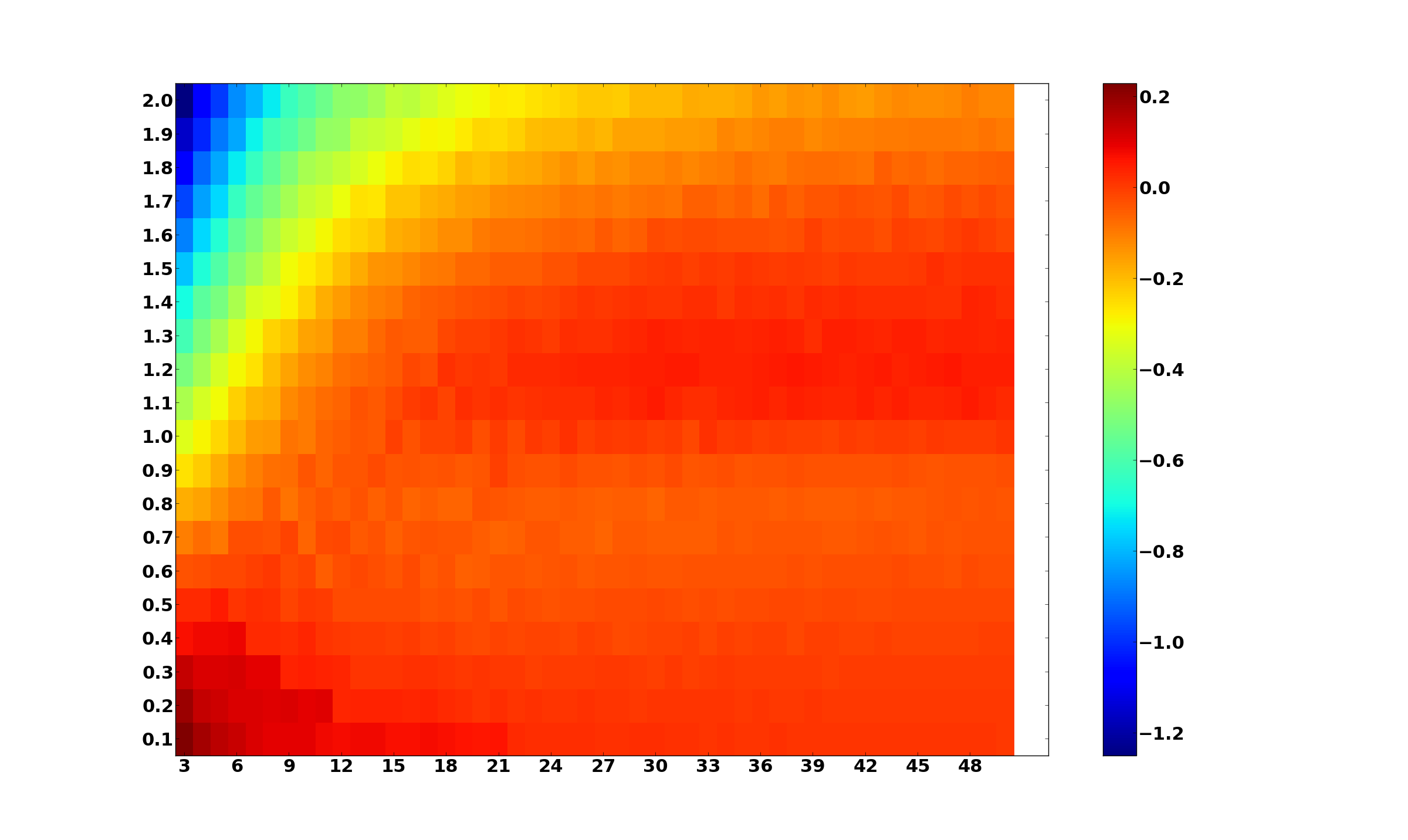}
        \end{subfigure}
        \\
        \begin{subfigure}[b]{0.4\textwidth}
            \centering
            \includegraphics[width=\textwidth]{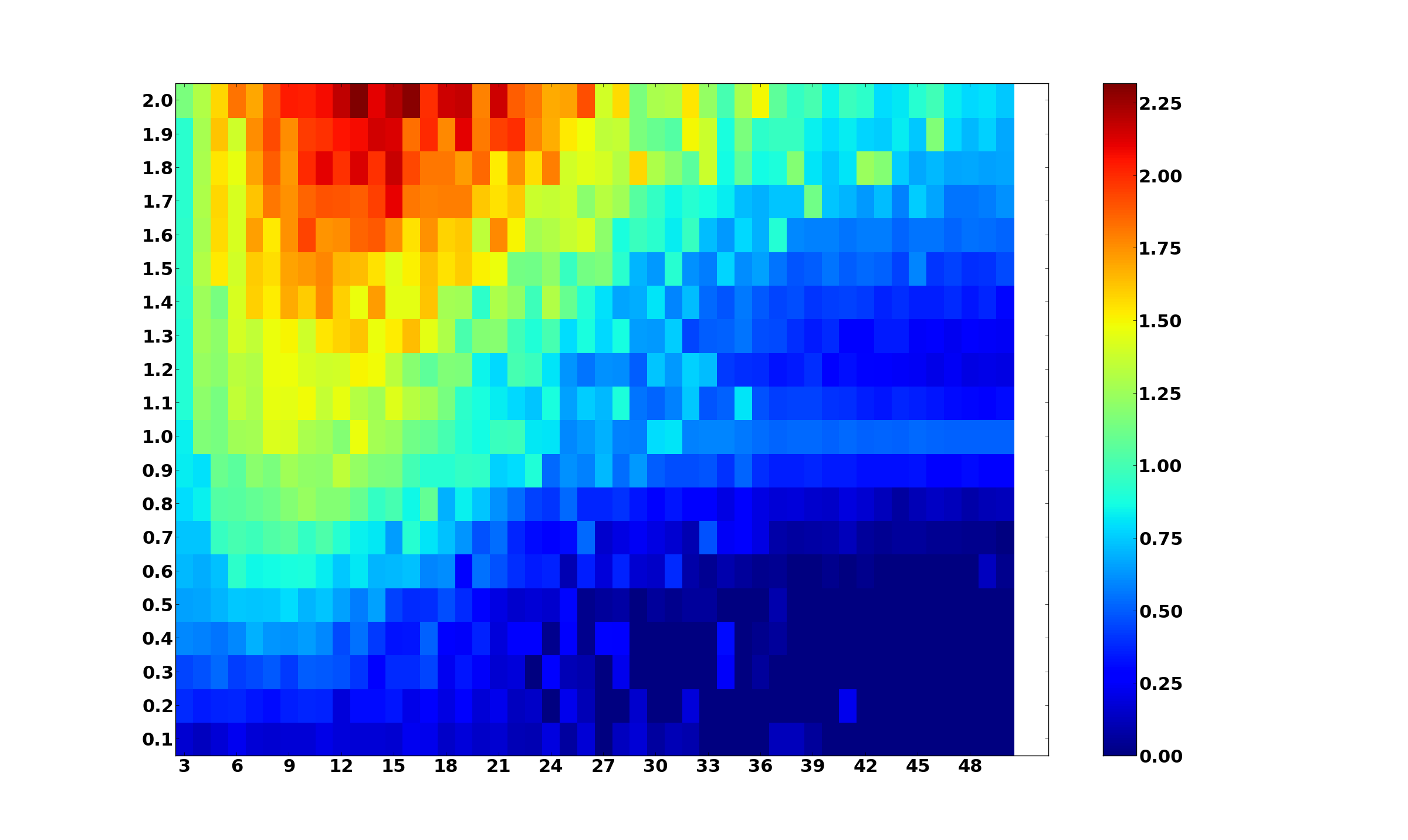}
            \caption{Counting}
        \end{subfigure}%
        ~ 
        \begin{subfigure}[b]{0.4\textwidth}
            \centering
            \includegraphics[width=\textwidth]{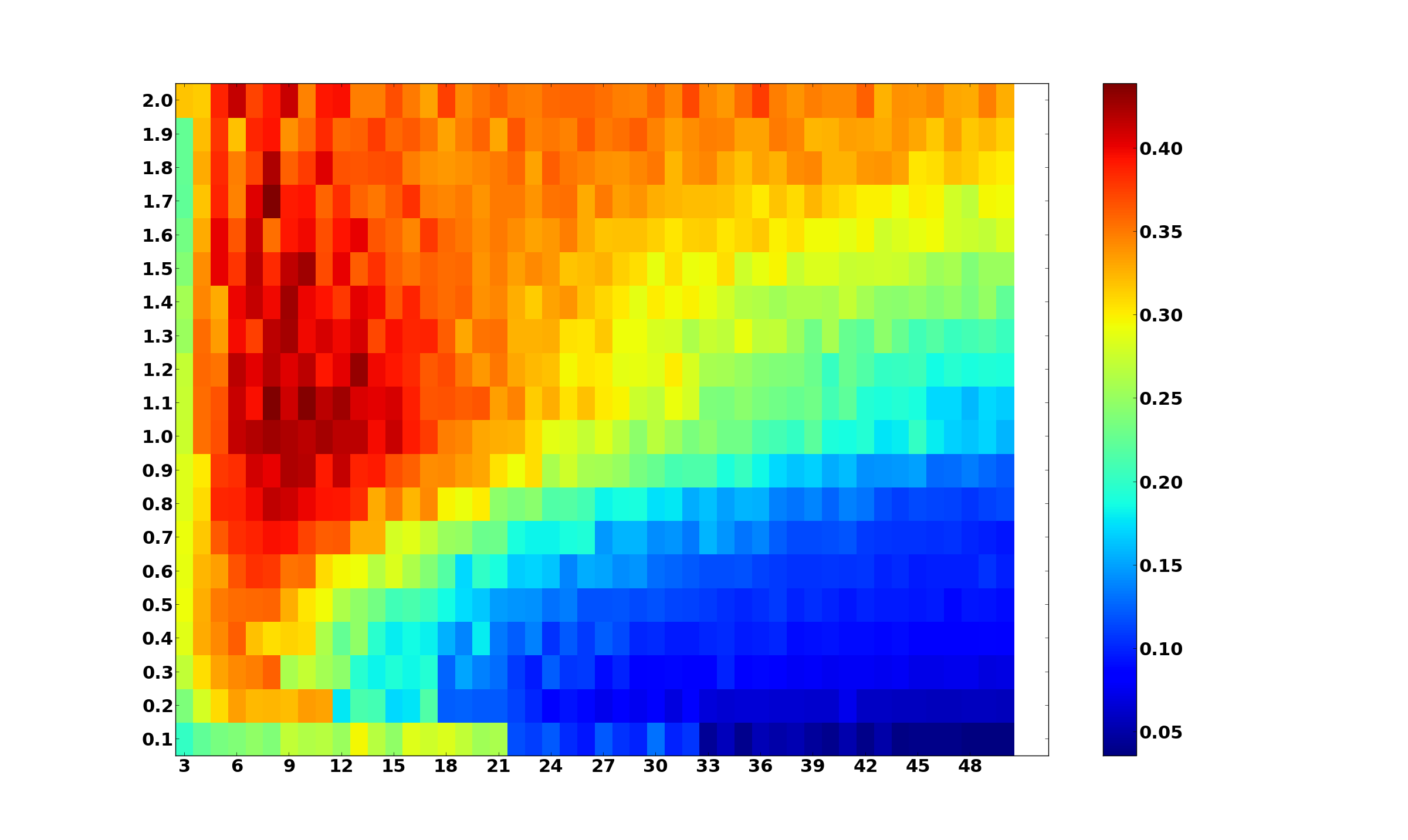}
            \caption{Bayes Mean}
        \end{subfigure}
        \begin{subfigure}[b]{0.4\textwidth}
            \centering
            \includegraphics[width=\textwidth]{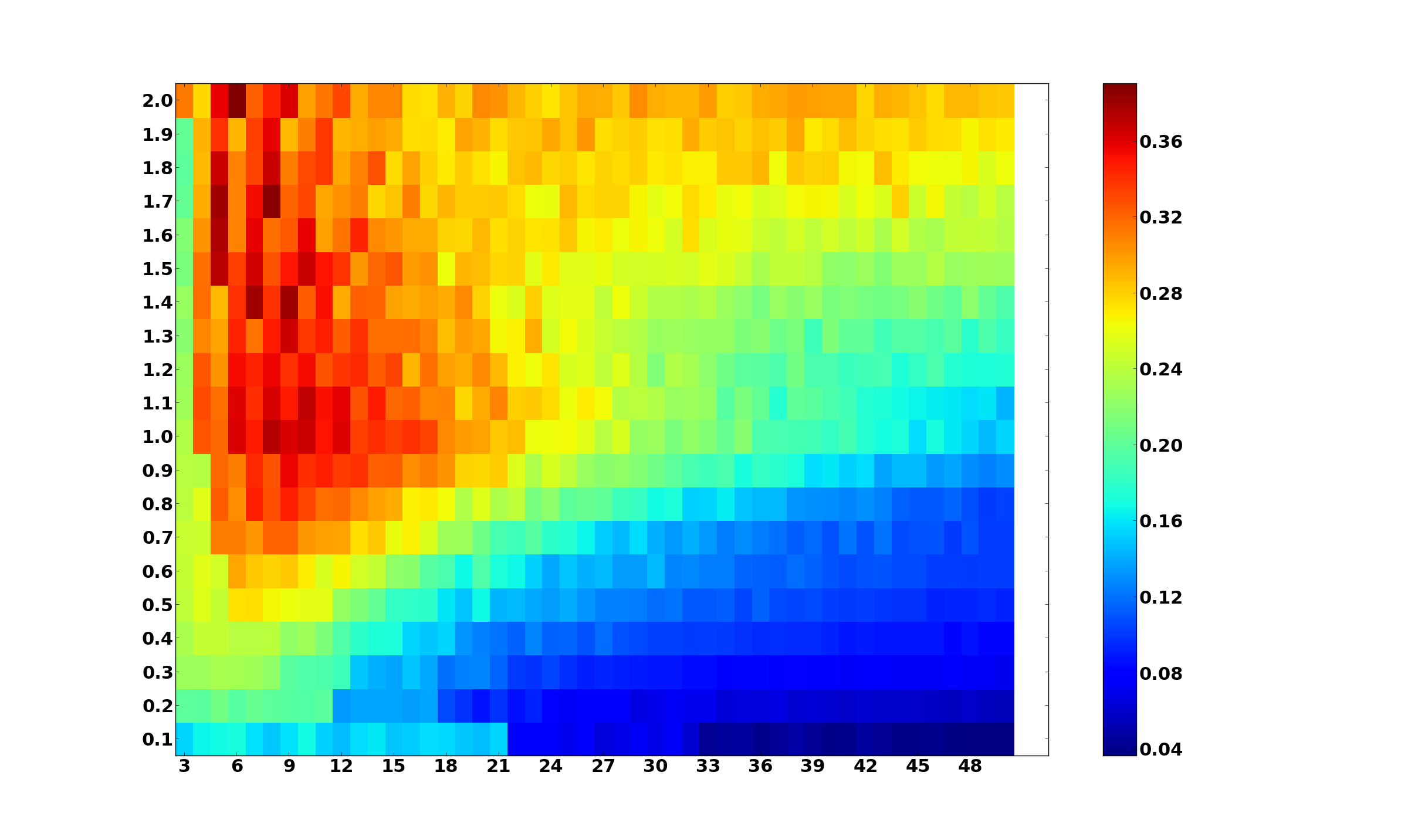}
            \caption{Bayes Mode}
        \end{subfigure}        
        \caption{Top Row: Color is given by the inferred values of r minus the true value for Counting, Bayesian Mean, and Bayesian Mode (left to right) for 1000 simulations at each pair $(N, r)$. Bottom: Standard deviations of inferred values for same methods. Initial point: $(a, N-a)$ where $a = \text{max}\left(1, \frac{N r}{r+1} \right)$. Although the deviation pattern looks similar, the magnitude is very different from counting to either inference estimate.}
        \label{moran_sims_balance}
\end{figure} 
\end{landscape}


\begin{landscape}
\centering
\begin{figure}[h]
        \begin{subfigure}[b]{0.4\textwidth}
            \centering
            \includegraphics[width=\textwidth]{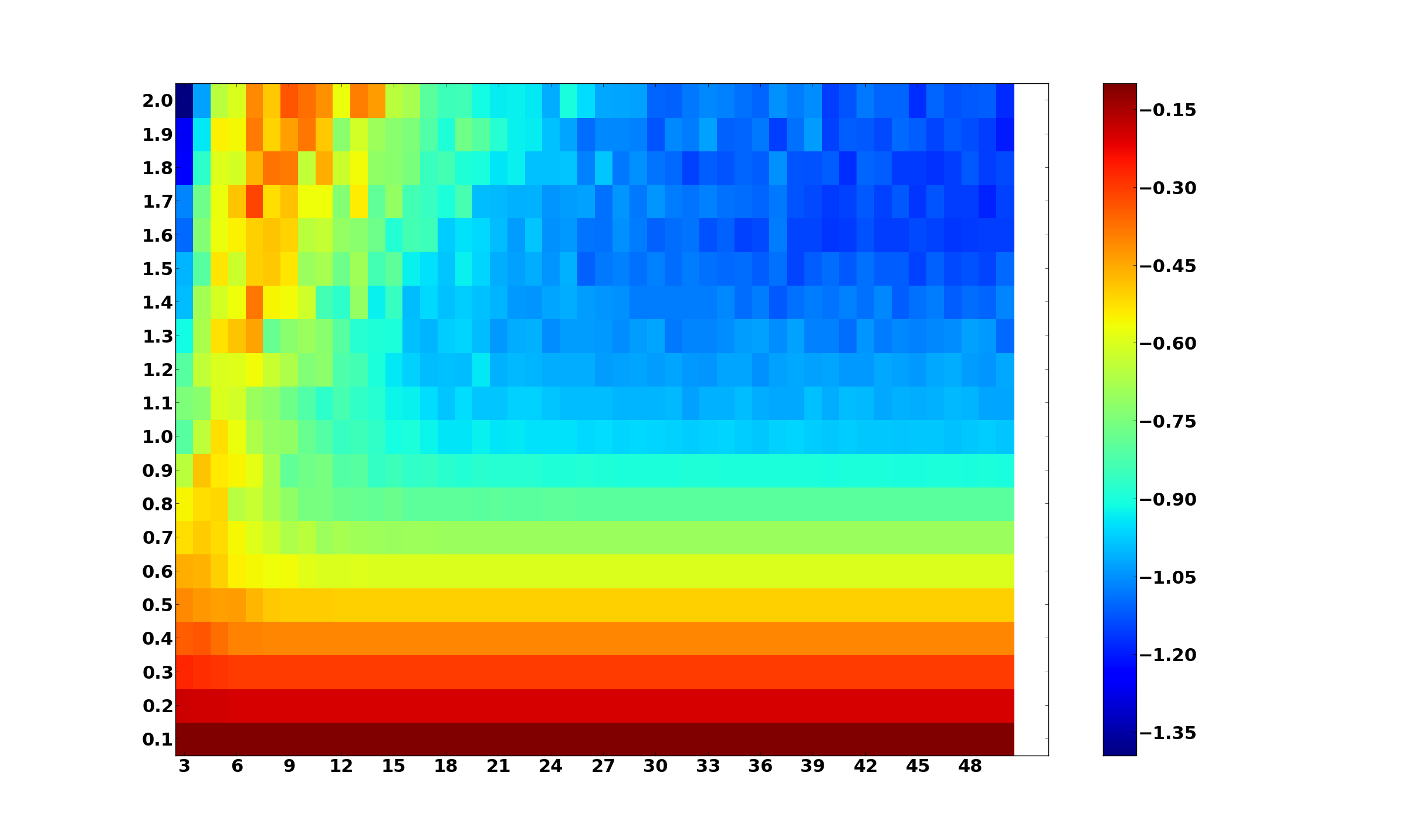}
        \end{subfigure}%
        ~ 
        \begin{subfigure}[b]{0.4\textwidth}
            \centering
            \includegraphics[width=\textwidth]{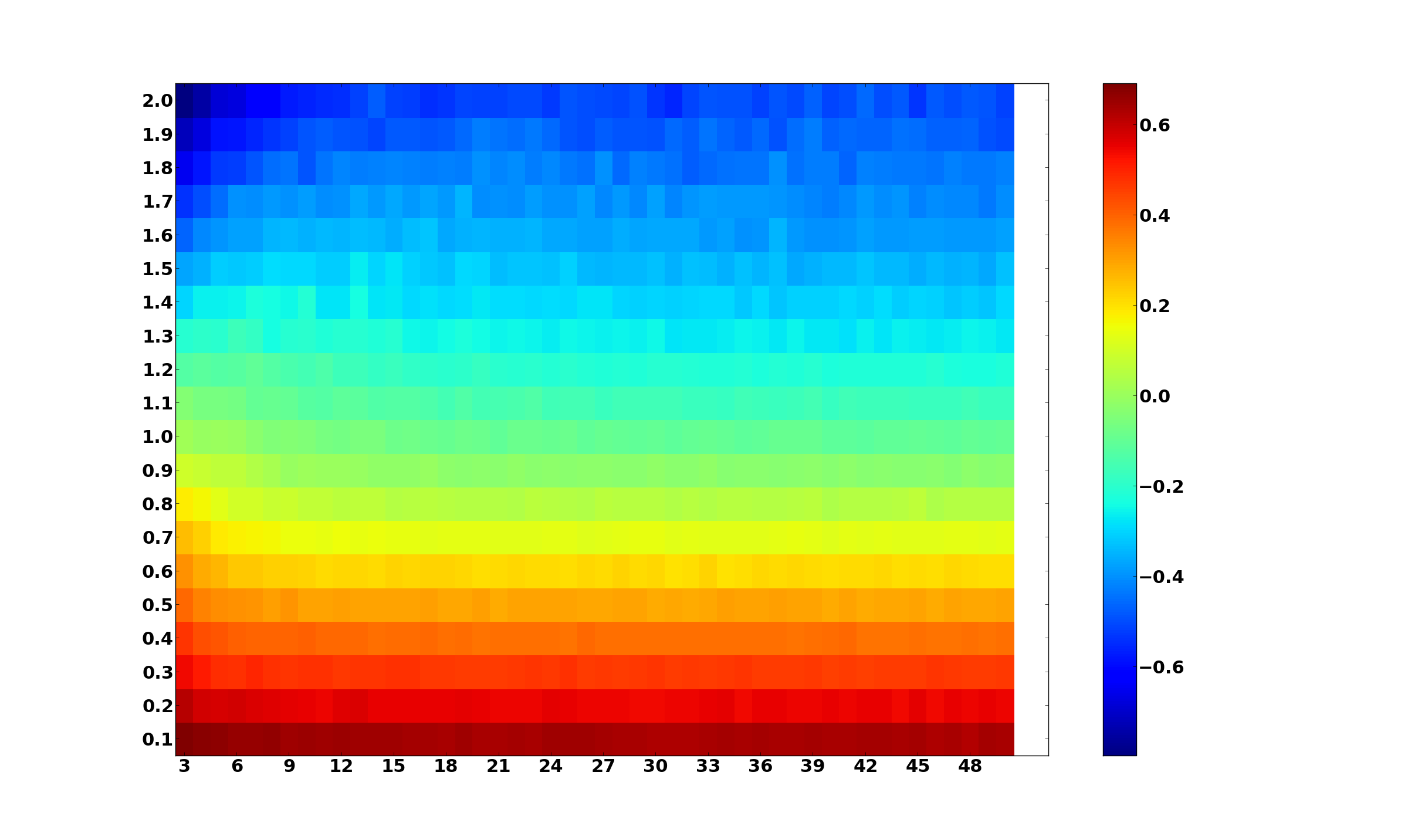}
        \end{subfigure}
        \begin{subfigure}[b]{0.4\textwidth}
            \centering
            \includegraphics[width=\textwidth]{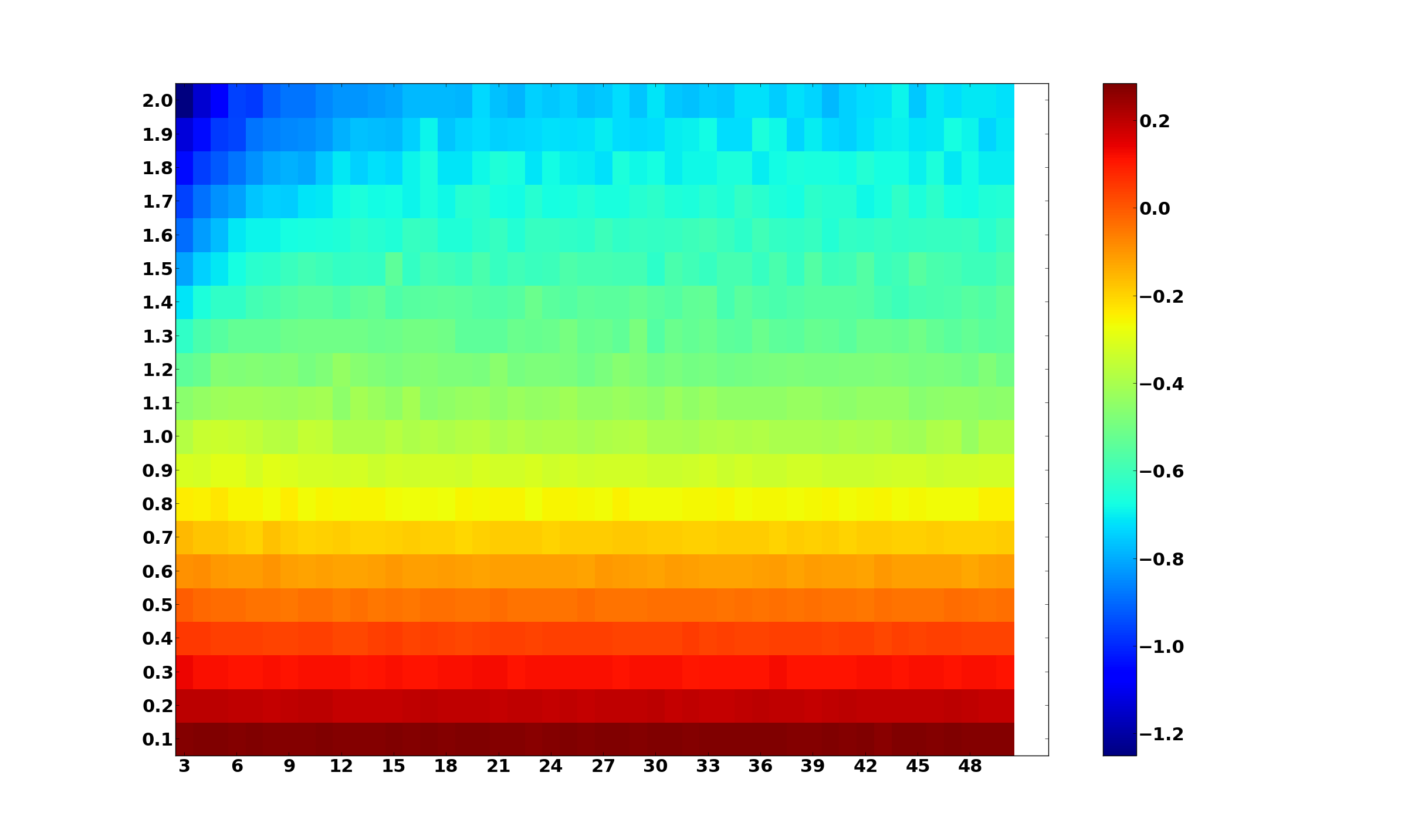}
        \end{subfigure}
        \\
        \begin{subfigure}[b]{0.4\textwidth}
            \centering
            \includegraphics[width=\textwidth]{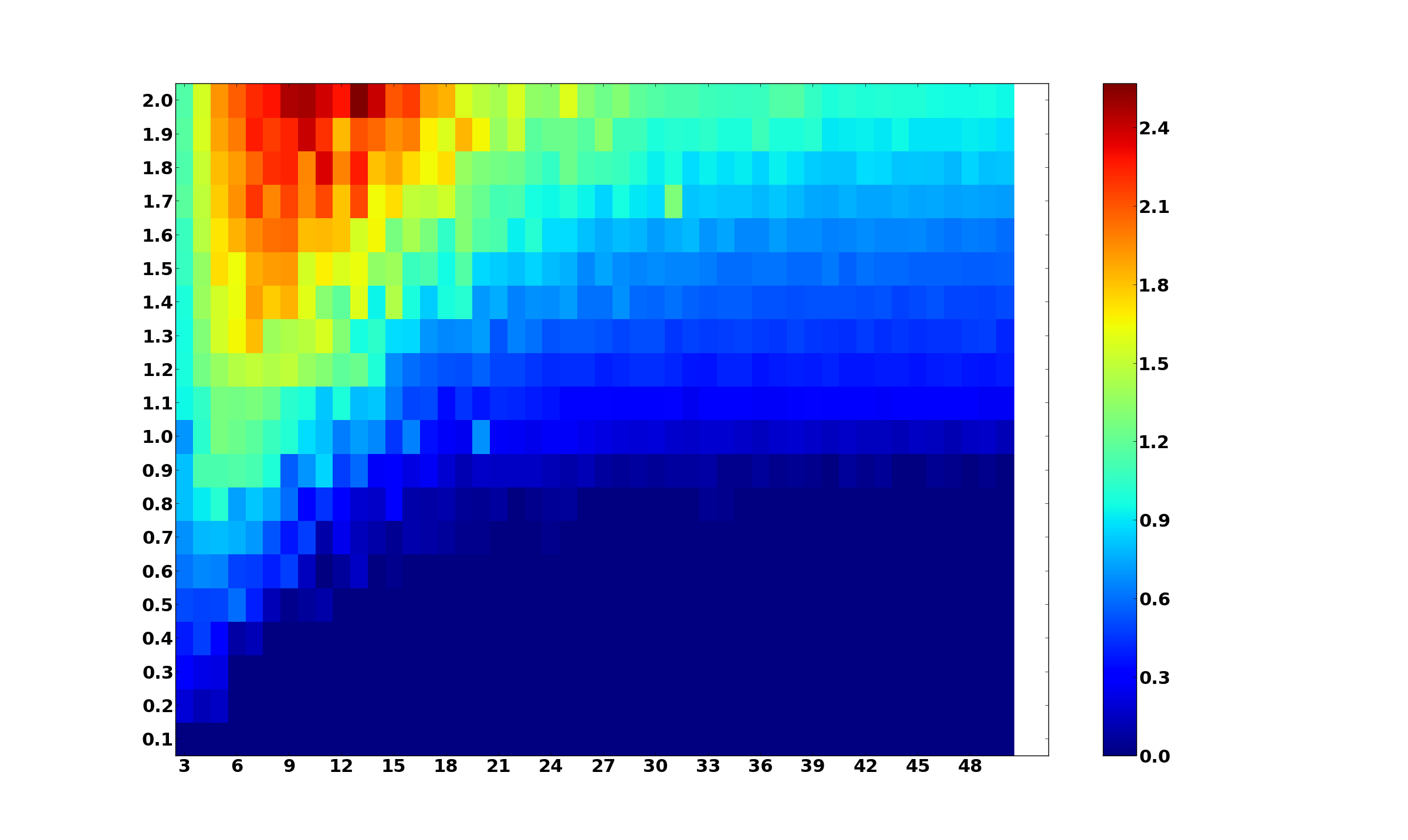}
            \caption{Counting}
        \end{subfigure}%
        ~ 
        \begin{subfigure}[b]{0.4\textwidth}
            \centering
            \includegraphics[width=\textwidth]{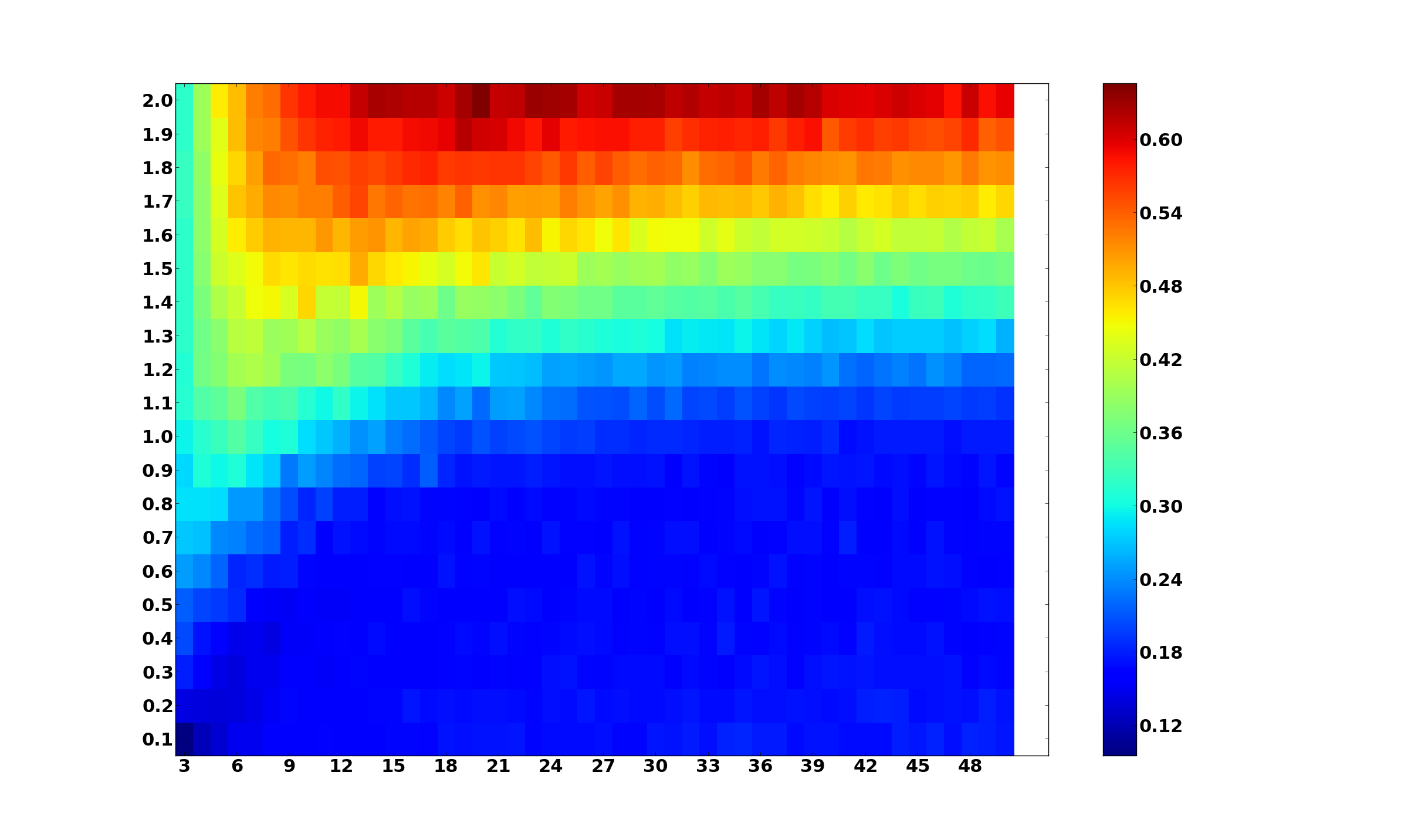}
            \caption{Bayes Mean}
        \end{subfigure}
        \begin{subfigure}[b]{0.4\textwidth}
            \centering
            \includegraphics[width=\textwidth]{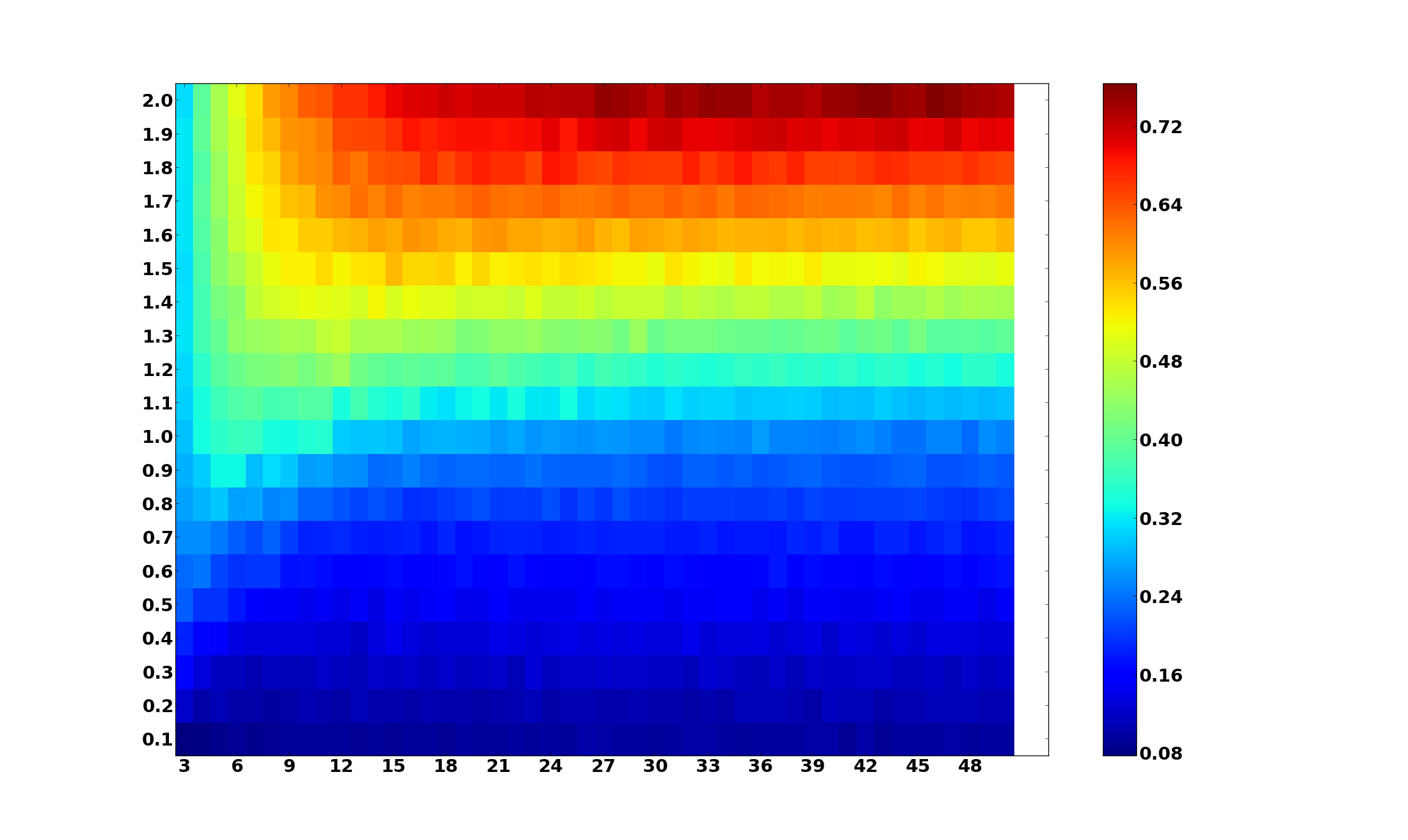}
            \caption{Bayes Mode}
        \end{subfigure}        
        \caption{Moran Process. Top Row: Color is given by the inferred values of r minus the true value for Counting, Bayesian Mean, and Bayesian Mode (left to right) for 1000 simulations at each pair $(N, r)$. Bottom: Standard deviations of inferred values for same methods. Initial point: $(1, N-1)$ (i.e. a single mutant). Although the deviation pattern looks similar, the magnitude is very different from counting to either inference estimate.}
        \label{moran_sims_single_mutant}
\end{figure} 
\end{landscape}

(1) Starting states adjacent to absorbing states: For all values of $r$ and $N$, starting states close to either absorbing state have a tendency to produce short trajectories, yielding little information for either method due to a shortage of data, especially if near a state that favors the true value of $r$. Moreover, for the counting method, estimates can be completely nonsensical since there may be zero births of one of the two types, leading to estimates of zero or infinity. Inference performs much better in this situation, both on such trajectories and generally, being more accurate and having a smaller variance, but depending heavily on the prior distribution.

(2) Small $N$: for small values of $N$, 
trajectory length can be short since the process can converge quickly. Nevertheless, both methods perform well if the starting state if is not adjacent to an absorbing state and $N>6$, with inference performing better generally. The transition probabilities for small populations are more skewed (as functions of $r$) than for large populations (e.g. compare $1/(1+r)$ with $1/(1+30r)$). This means that observations for small populations can have a larger impact on the inferred value than for larger populations, and that more observations are required to lower variance. Figure \ref{moran_sims_single_mutant} shows that smaller population sizes have greater average deviation from the true value of $r$ and larger variances in inferred value.

(3) $r$ close to 1: For values of $r$ close to 1, the lower variance of inference produces significantly better estimates. Because trajectories can ultimately converge to either absorbing state with relatively similar probabilities (depending on the starting state), the additional information obtained by inference in each transition produces far better estimates. For $r=1$, the variance from counting can be ten times greater than the variance for inference in addition to inference producing a more accurate mean value. Note also that for values of $r >> 1$, a prior that has more weight away from $r=1$ can improve estimates, and the choice of prior may favor this case.

Although it appears that counting performs similarly to inference in the extreme case of a single mutant, it is critical to note that these values are only computed in cases that counting could give an estimate. That is, in many cases, especially for extreme relative fitnesses and starting states, only one type registers any replication events, which can lead to estimates of infinite relative fitness. This affect is slightly less extreme for processes that separate death and birth in accounting.

\subsection{Separating Birth and Death}

Counting suffers from the fact that birth and death events are only clear in the case that the population changes state. Since the most likely transition is very often to stay in the current state, even if the information gained from such a transition is small, inference has a big advantage due to the sheer number of such transitions, especially when close to an absorbing state. Trajectories of the Moran process tend to oscillate near the absorbing state because while the higher fit type typically dominates the population, it is also most likely to be randomly selected for death. For instance, suppose $r = 3$ and $N=20$. Then the population will have a tendency to cycle in the states $(1, 19) \to (1, 20) \to (1, 19)$. The reason is that while and individual of type $B$ is far more likely to be chosen to reproduce because $r > 1$ and nearly all individuals are of type $B$, it is also unlikely that the lone individual of type $A$ will be chosen for death: $Pr(\text{chose $B$ individual for removal}) = 1 / (N+1)$. This can lead to a large number of $B$ replication events near fixation which can lead the estimate to drift upward near the end of the trajectory.

\begin{figure}[h]
    \centering
    \includegraphics[scale=0.28]{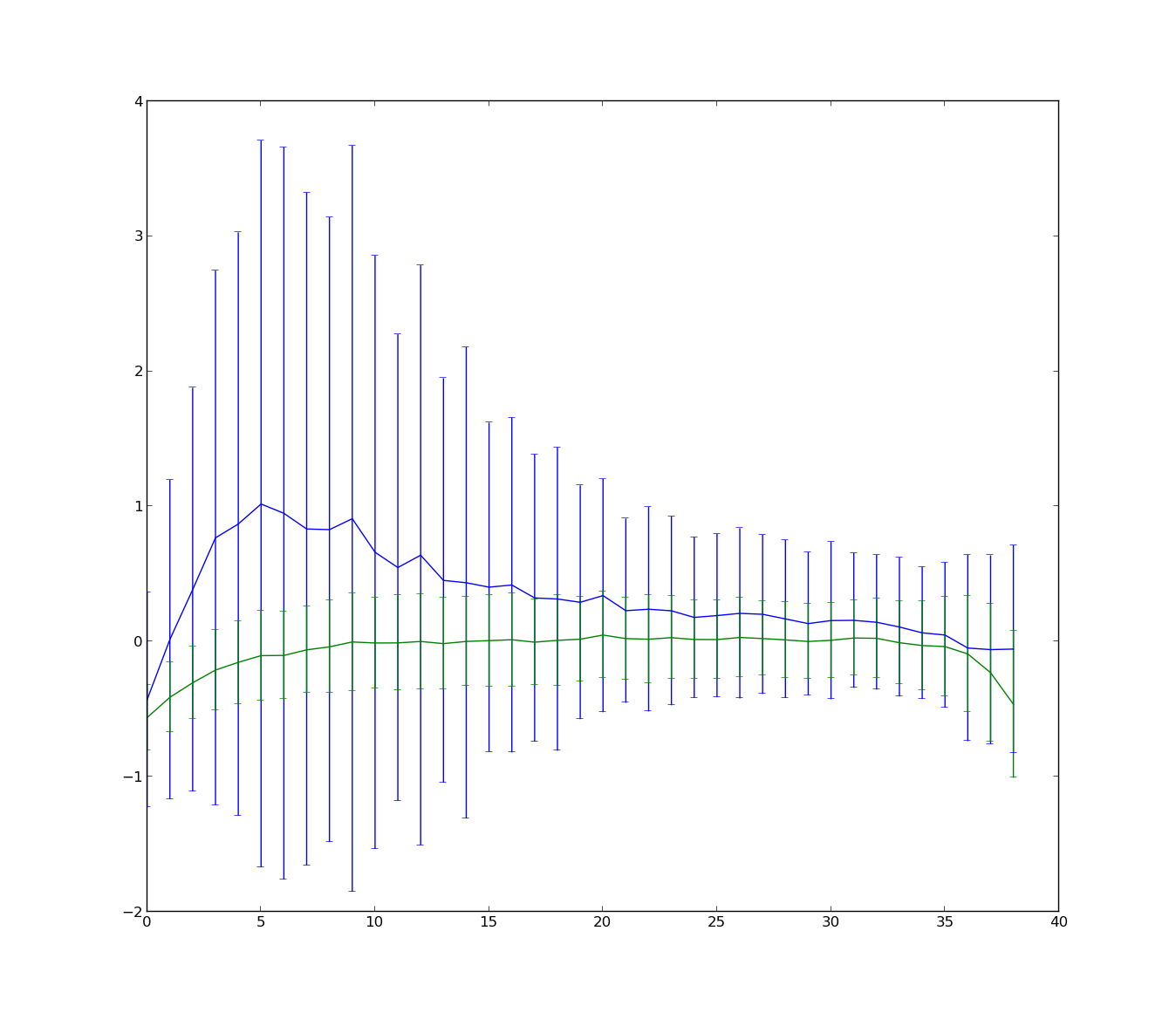}
    \caption{Comparison of Inverse counting (blue) and Bayesian mean (green). Curves indicate the deviation from the true value of $r$ and the standard deviation over 1000 simulations for $r = 1.9$ and population size $N$ along the x-axis.}
    \label{inverse_counting_bayesian_mean}
\end{figure}

To understand the strength of the effect of separating birth and death on the ability to infer fitness, compare the results from the Moran process to the modified process in which birth and death events are recorded separately. The method of counting also benefits from this process, but care must be taken. Naively proceeding as before is unwise because of the often quite small transition probabilities near absorbing states. To manage this over-abundance of births of the more fit type, we count weighted by population size. If an $A$ individual reproduces in a population in state $(a,b)$, this is counted as $1/a$; similarly $B$ births are counted as $1/b$. This is justified by the estimate given previously and produces much improved estimates (since the transitions in which the population stays in the same state now count as replication events and the large numbers of oscillator transitions near the absorbing states would lead to substantial overestimations). For the modified process, both methods perform well, significantly improving versus the Moran process. Inference still maintains a smaller variance, and overall has lower variance and more accurate estimates for this process than the Moran process. In general, counting performs poorly (see Figures \ref{inverse_counting_bayesian_mean} and \ref{fps_sims_balanced}).

\begin{landscape}
\centering
\begin{figure}[h]
        \begin{subfigure}[b]{0.4\textwidth}
            \centering
            \includegraphics[width=\textwidth]{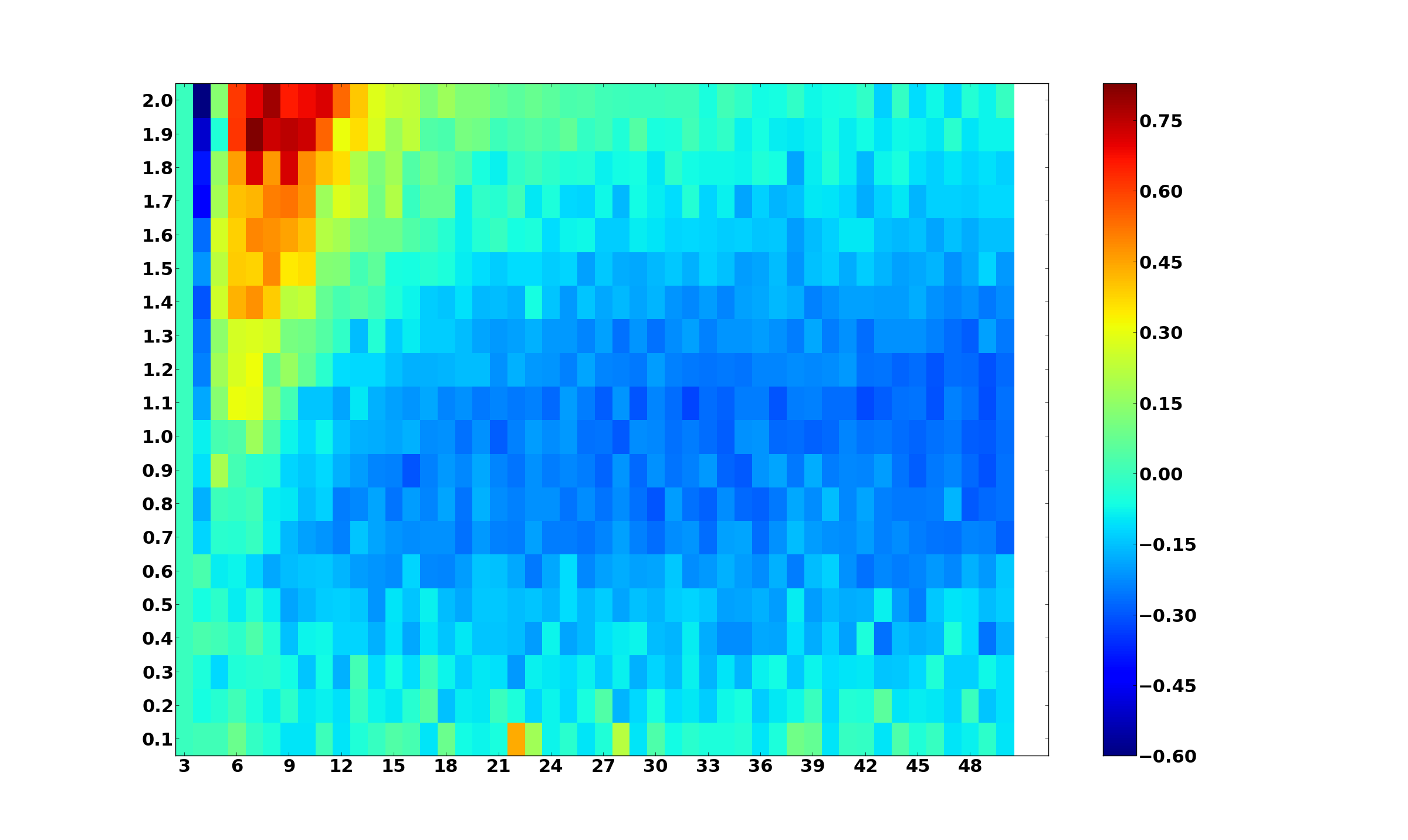}
        \end{subfigure}%
        ~ 
        \begin{subfigure}[b]{0.4\textwidth}
            \centering
            \includegraphics[width=\textwidth]{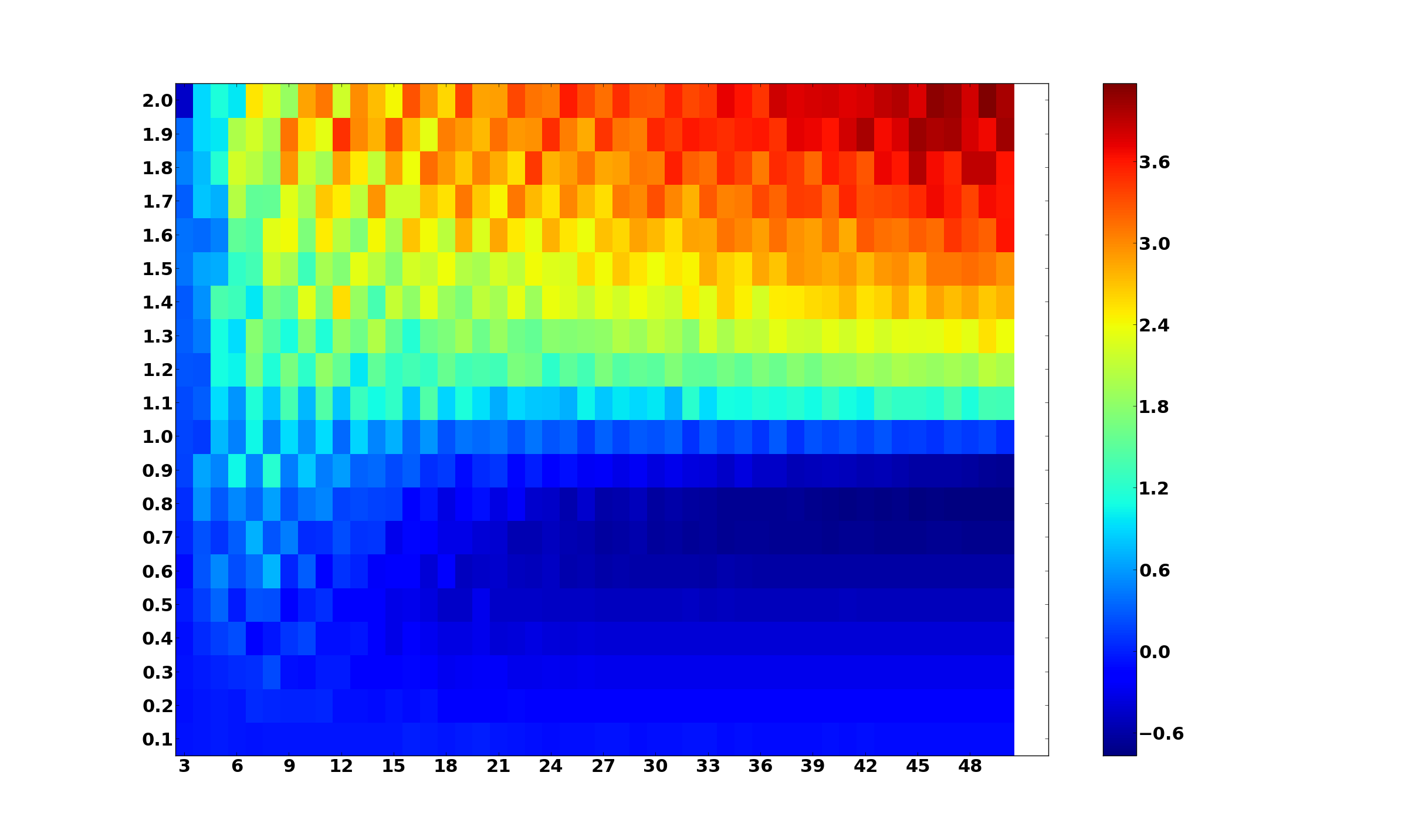}
        \end{subfigure}
        \begin{subfigure}[b]{0.4\textwidth}
            \centering
            \includegraphics[width=\textwidth]{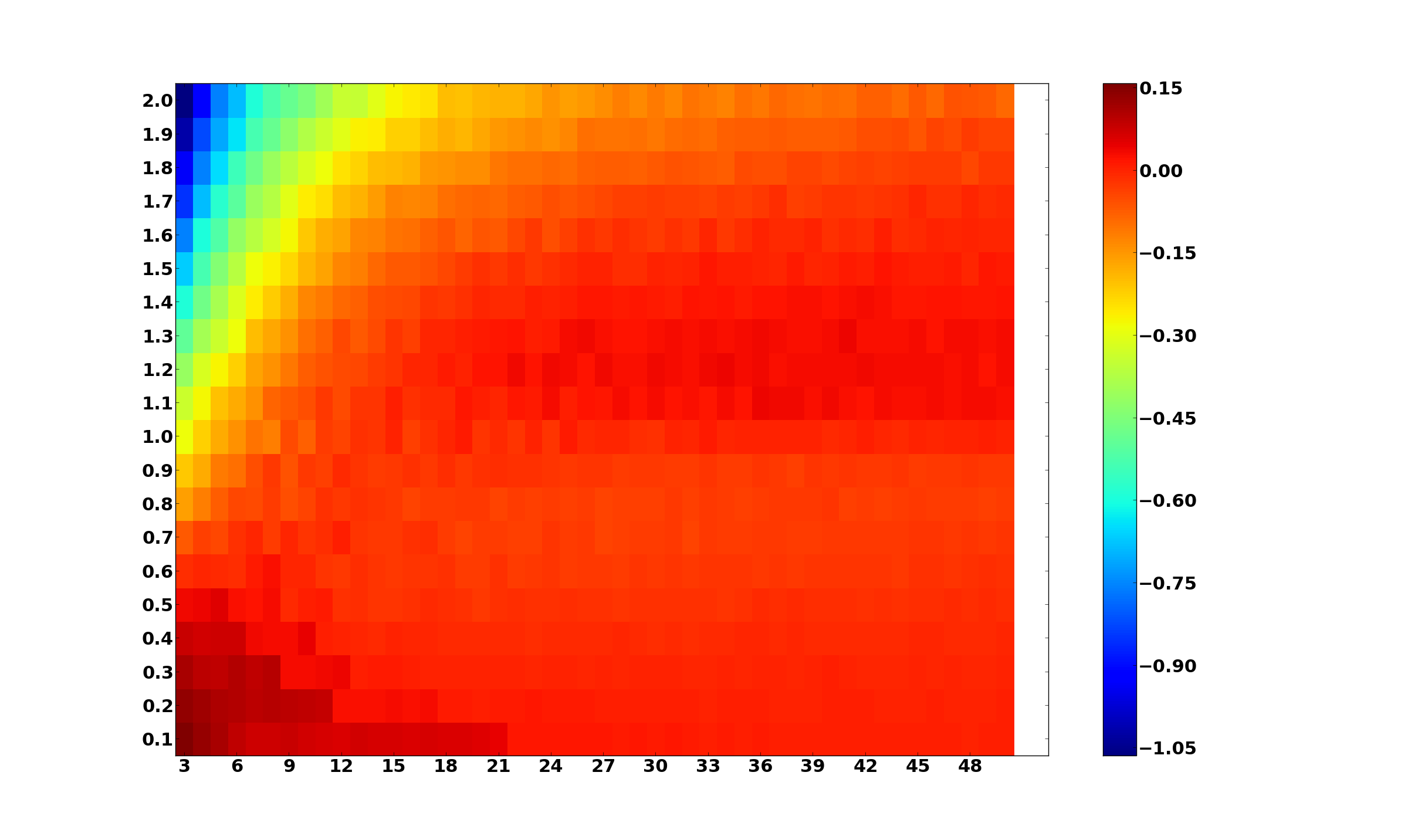}
        \end{subfigure}
        \\
        \begin{subfigure}[b]{0.4\textwidth}
            \centering
            \includegraphics[width=\textwidth]{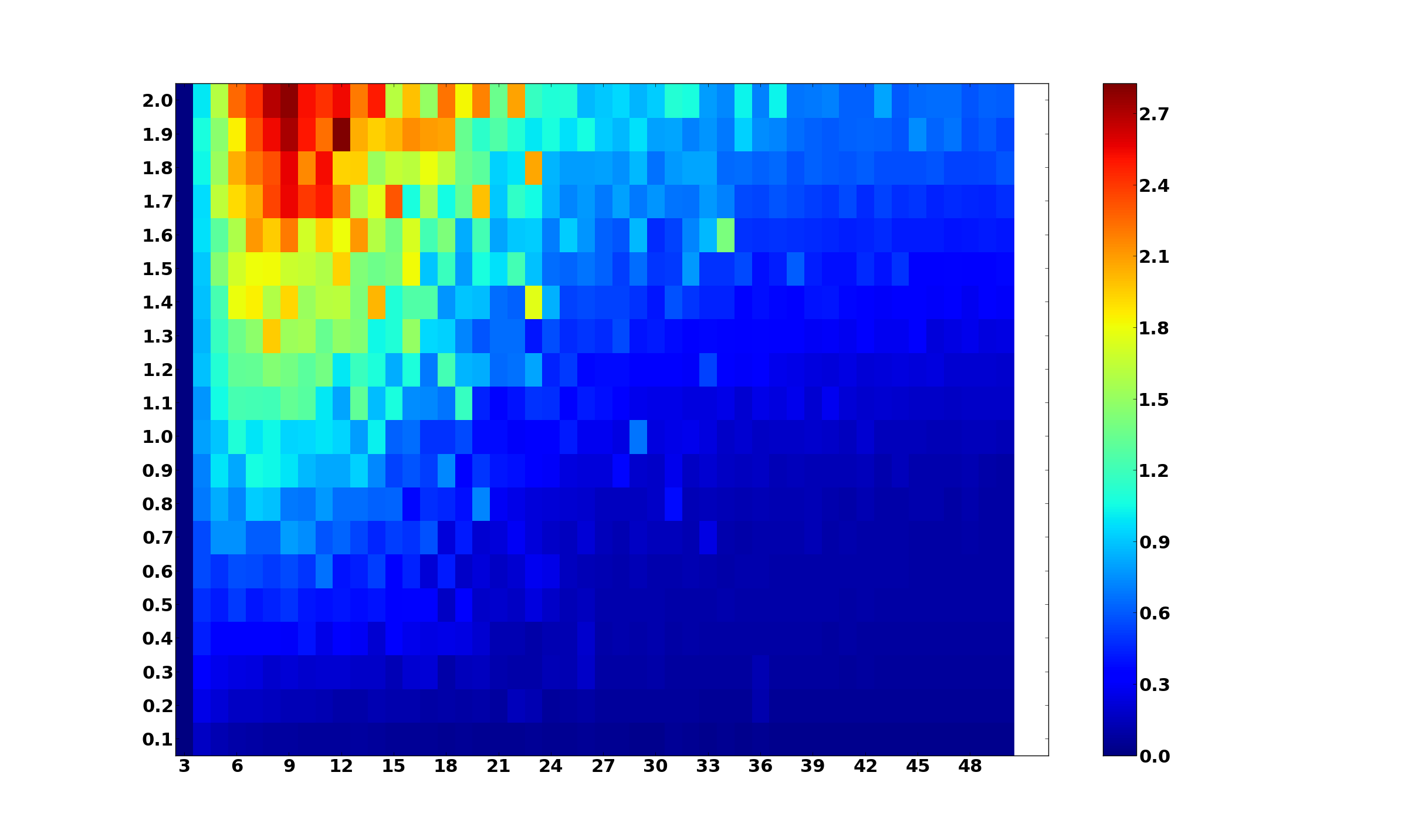}
            \caption{Inv. Counting}
        \end{subfigure}%
        ~ 
        \begin{subfigure}[b]{0.4\textwidth}
            \centering
            \includegraphics[width=\textwidth]{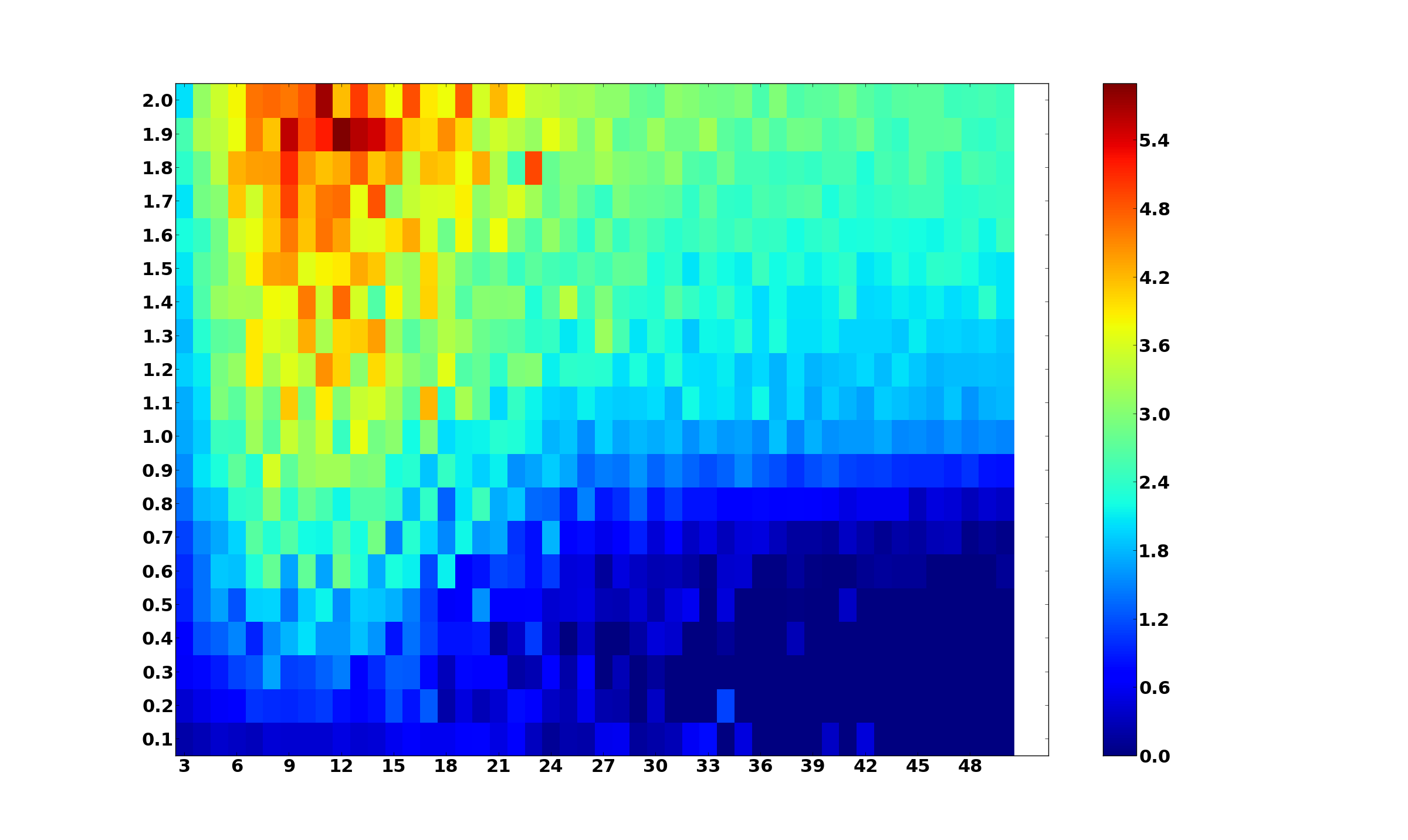}
            \caption{Counting}
        \end{subfigure}
        \begin{subfigure}[b]{0.4\textwidth}
            \centering
            \includegraphics[width=\textwidth]{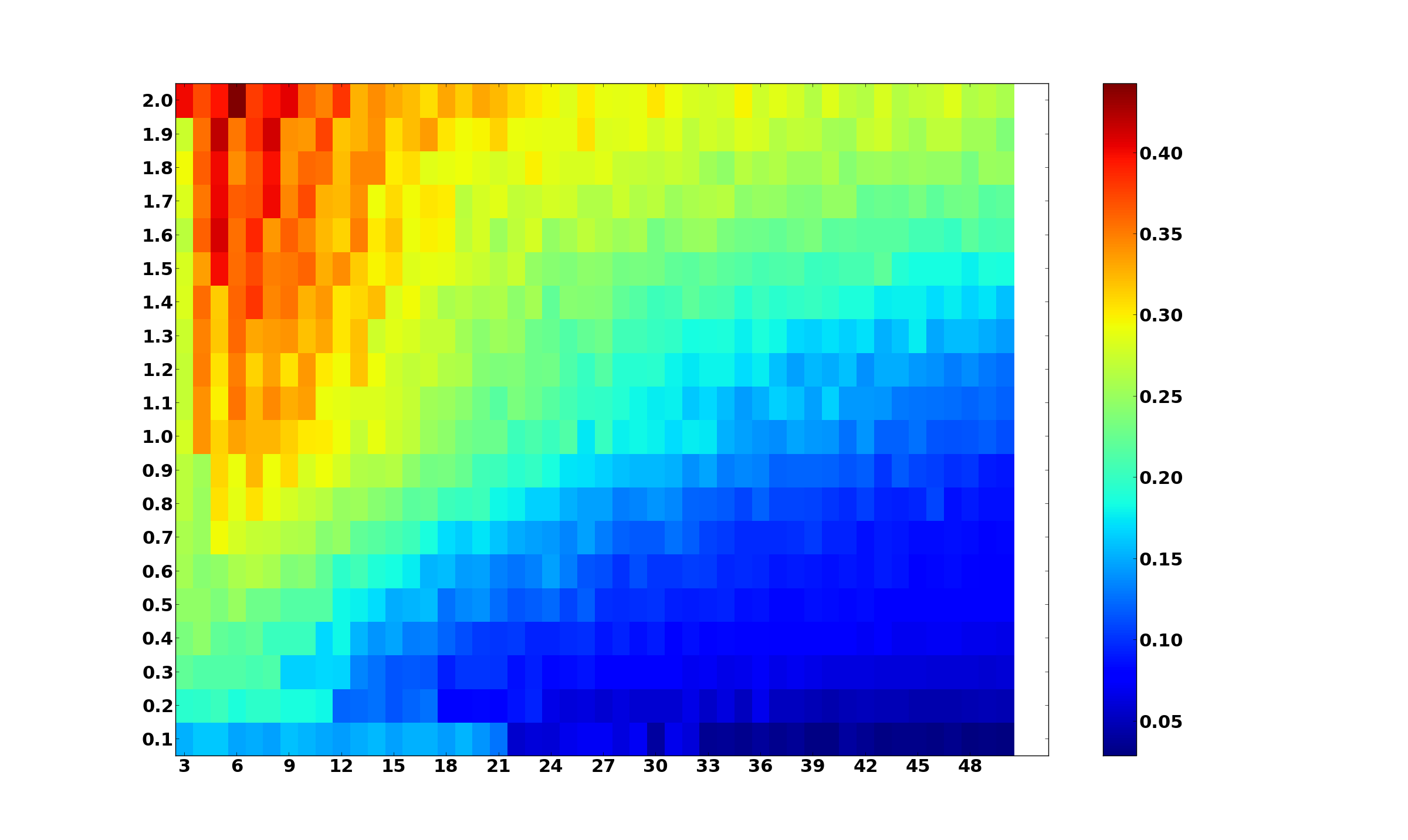}
            \caption{Bayes Mode}
        \end{subfigure}        
        \caption{Separated Moran Process. Top Row: Color is given by the inferred values of r minus the true value for Inverse Counting, Counting, and Bayesian Mode (left to right) for 1000 simulations at each pair $(N, r)$. Bottom: Standard deviations of inferred values for same methods. Initial point: $(a, N-a)$ where $a = \text{max}\left(1, \frac{N r}{r+1} \right)$. Although the deviation pattern looks similar, the magnitude is very different from counting to either inference estimate.}
        \label{fps_sims_balanced}
\end{figure} 
\end{landscape}

\subsection{Sampling}

As a practical matter, the full trajectory of may not be available or practically obtainable, so we also compare the two methods by sampling the trajectories for a subset of transitions in each trajectory. For sample sizes of 10 and 20 (with the full trajectory used if shorter than the sample size), the inference method suffers relatively little and mostly in variance in comparison to counting, which performs poorly even in cases where $N$ is large and the starting state is near the central state. So as a practical means of estimating fitness, it is not necessary to know the population trajectory to a high degree of precision, nor is the method dependent on the fact that there is dependence between states and observations. In other words, the fact that an observation of a transition from $(a,b)$ to $(a+1, b-1)$ is dependent on having had observations that got the population to the state $(a,b)$ in the first place would lead nonzero indices of $\alpha$ and $\beta$ to have nonzero neighbors (but this is not the case for a sampled trajectory).

\begin{landscape}
\centering
\begin{figure}[h]
        \begin{subfigure}[b]{0.4\textwidth}
            \centering
            \includegraphics[width=\textwidth]{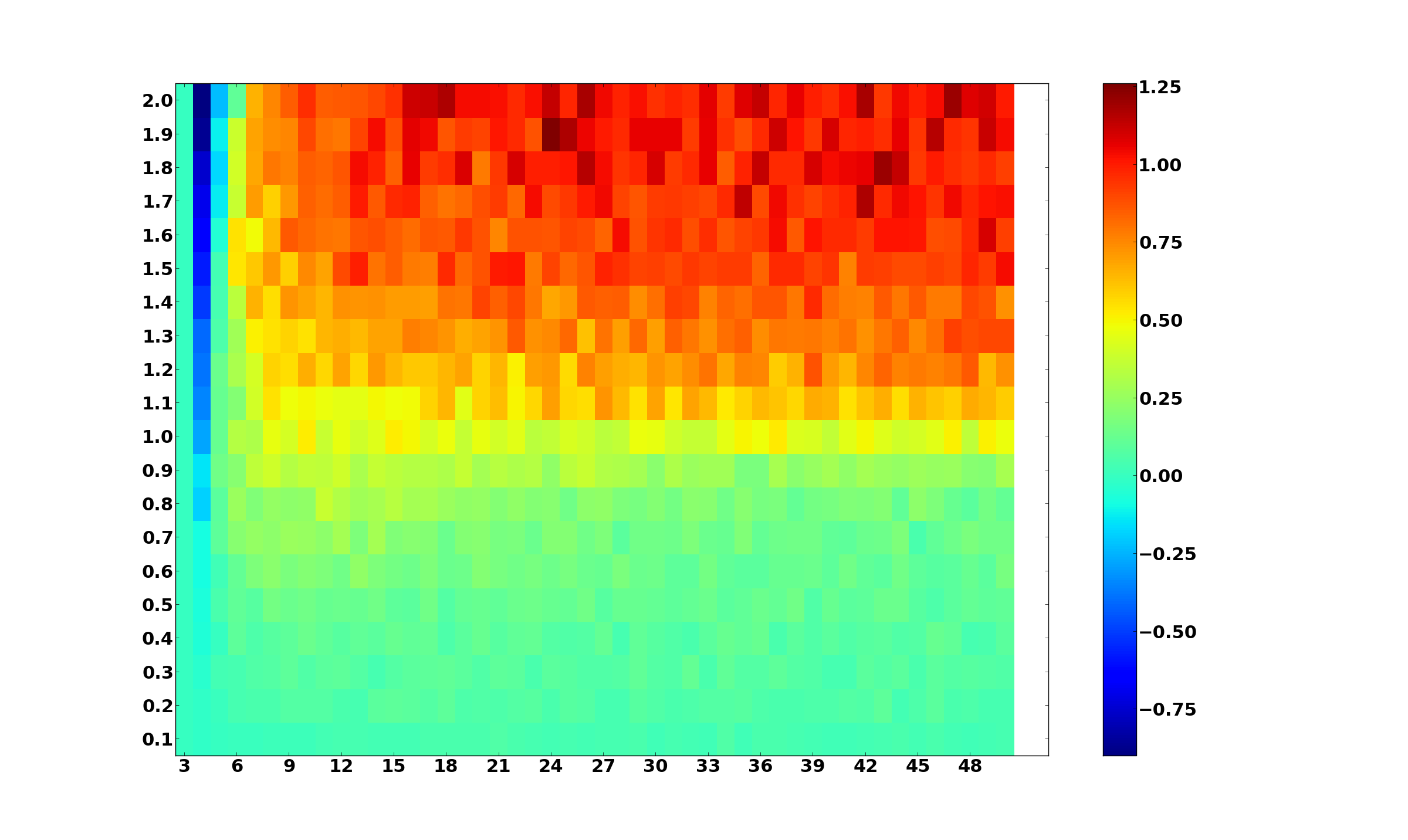}
        \end{subfigure}%
        ~ 
        \begin{subfigure}[b]{0.4\textwidth}
            \centering
            \includegraphics[width=\textwidth]{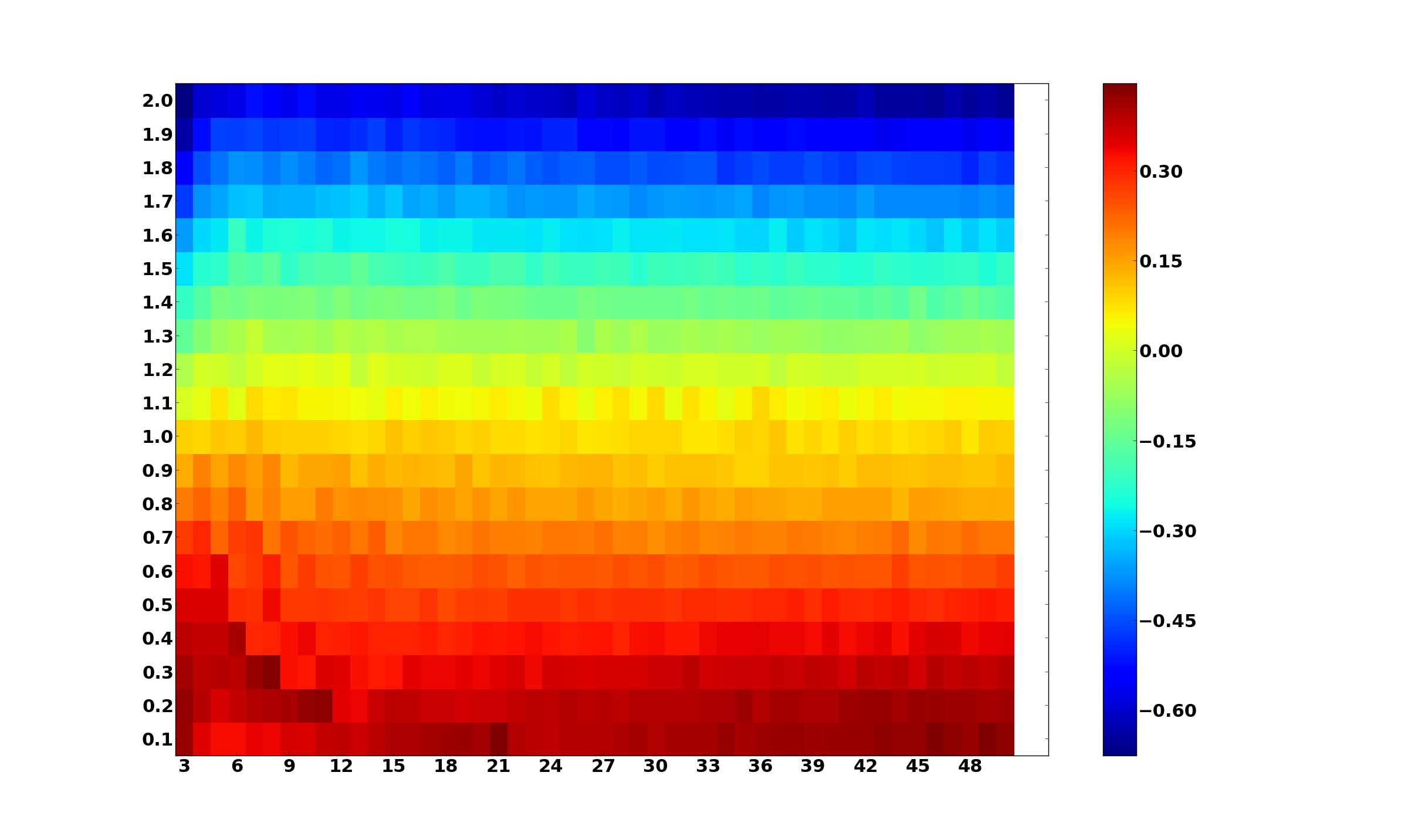}
        \end{subfigure}
        \begin{subfigure}[b]{0.4\textwidth}
            \centering
            \includegraphics[width=\textwidth]{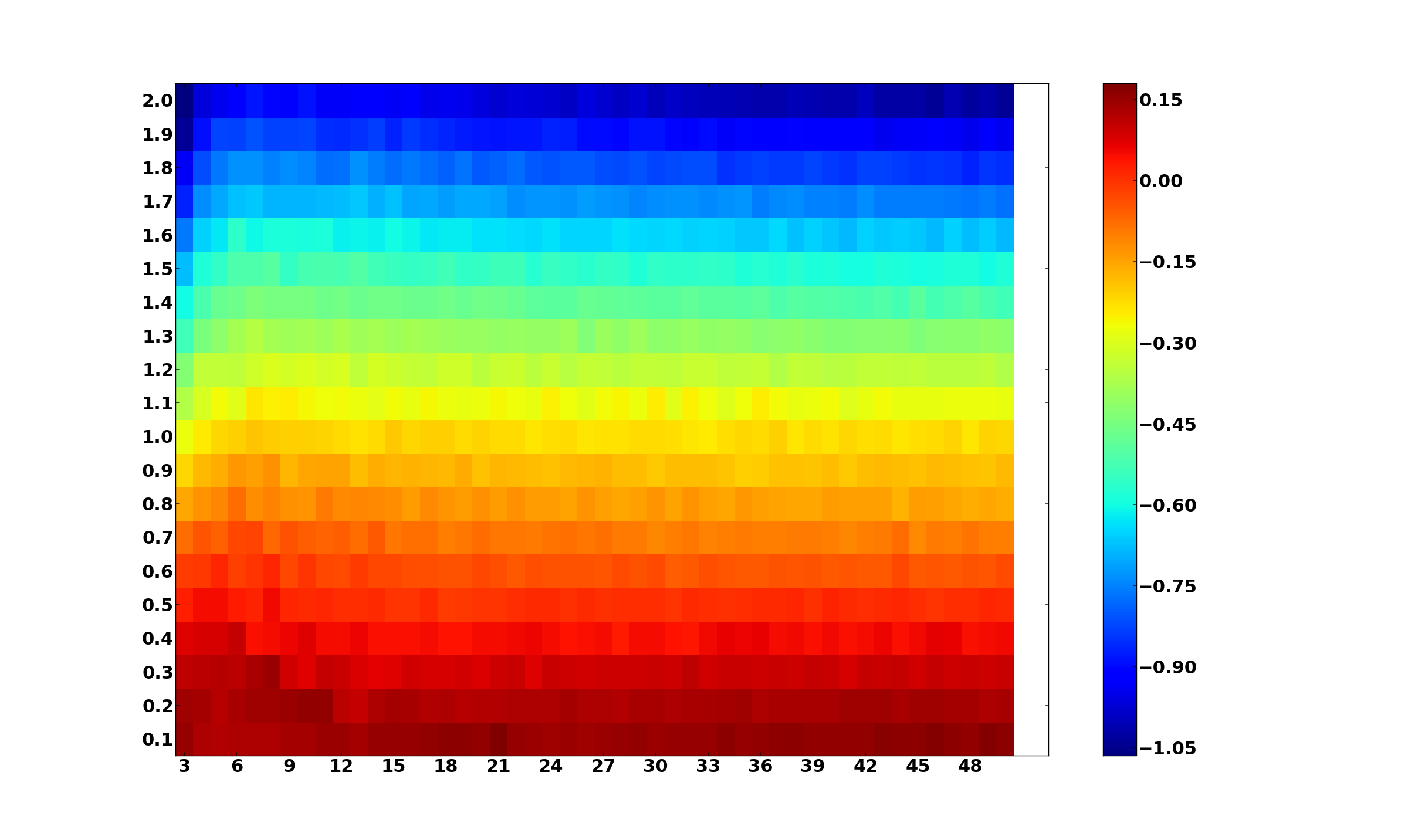}
        \end{subfigure}
        \\
        \begin{subfigure}[b]{0.4\textwidth}
            \centering
            \includegraphics[width=\textwidth]{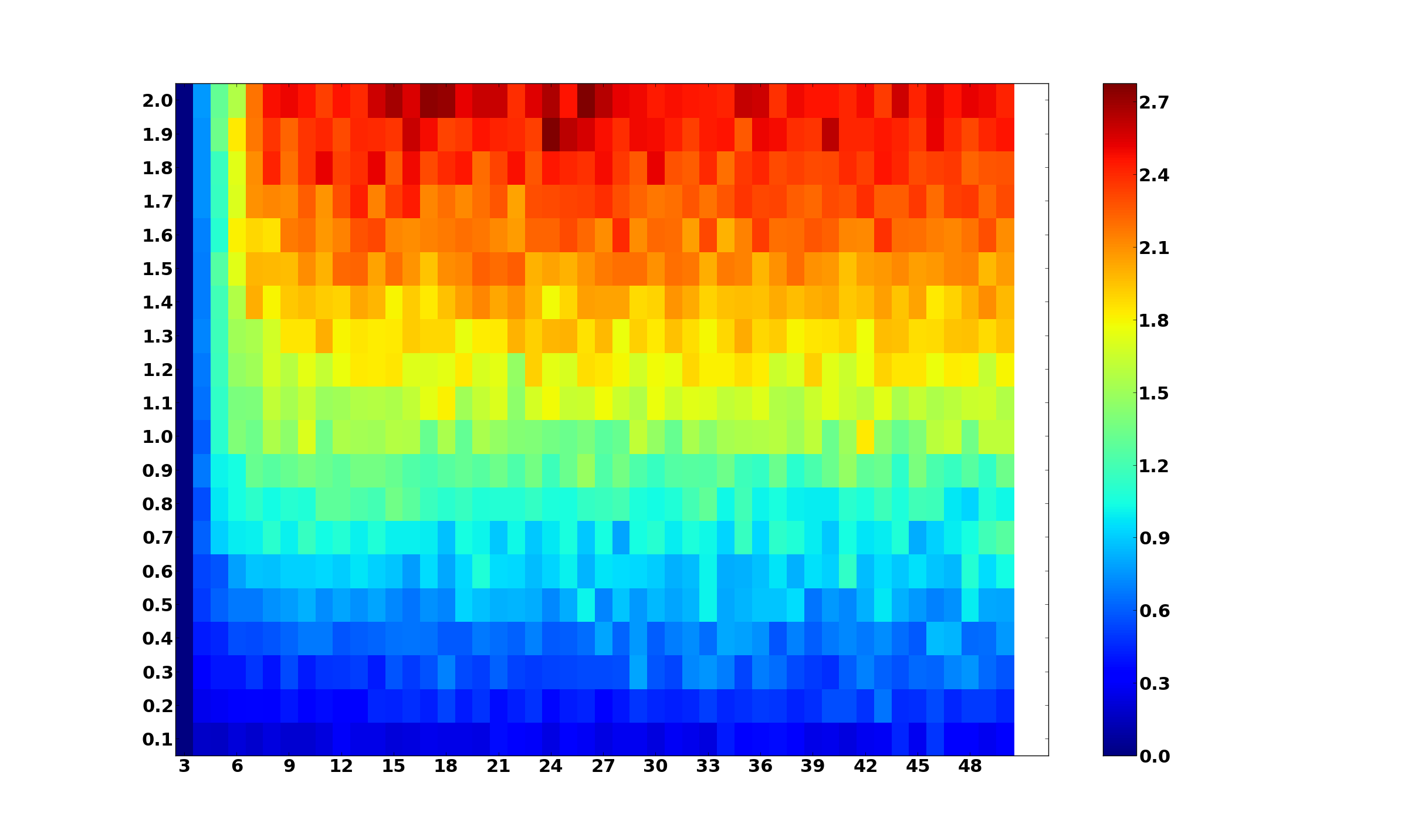}
            \caption{Inv. Counting}
        \end{subfigure}%
        ~ 
        \begin{subfigure}[b]{0.4\textwidth}
            \centering
            \includegraphics[width=\textwidth]{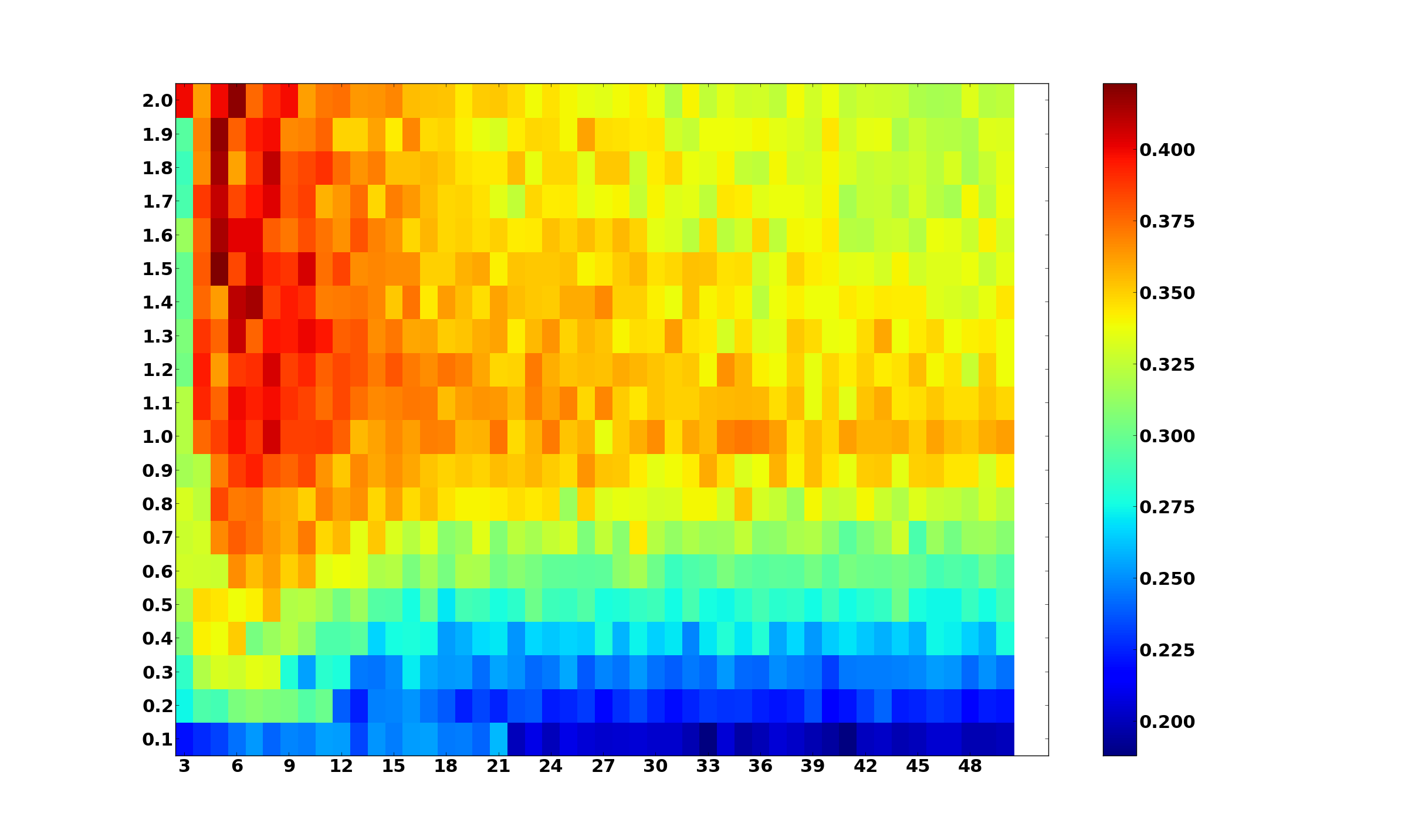}
            \caption{Bayes Mean}
        \end{subfigure}
        \begin{subfigure}[b]{0.4\textwidth}
            \centering
            \includegraphics[width=\textwidth]{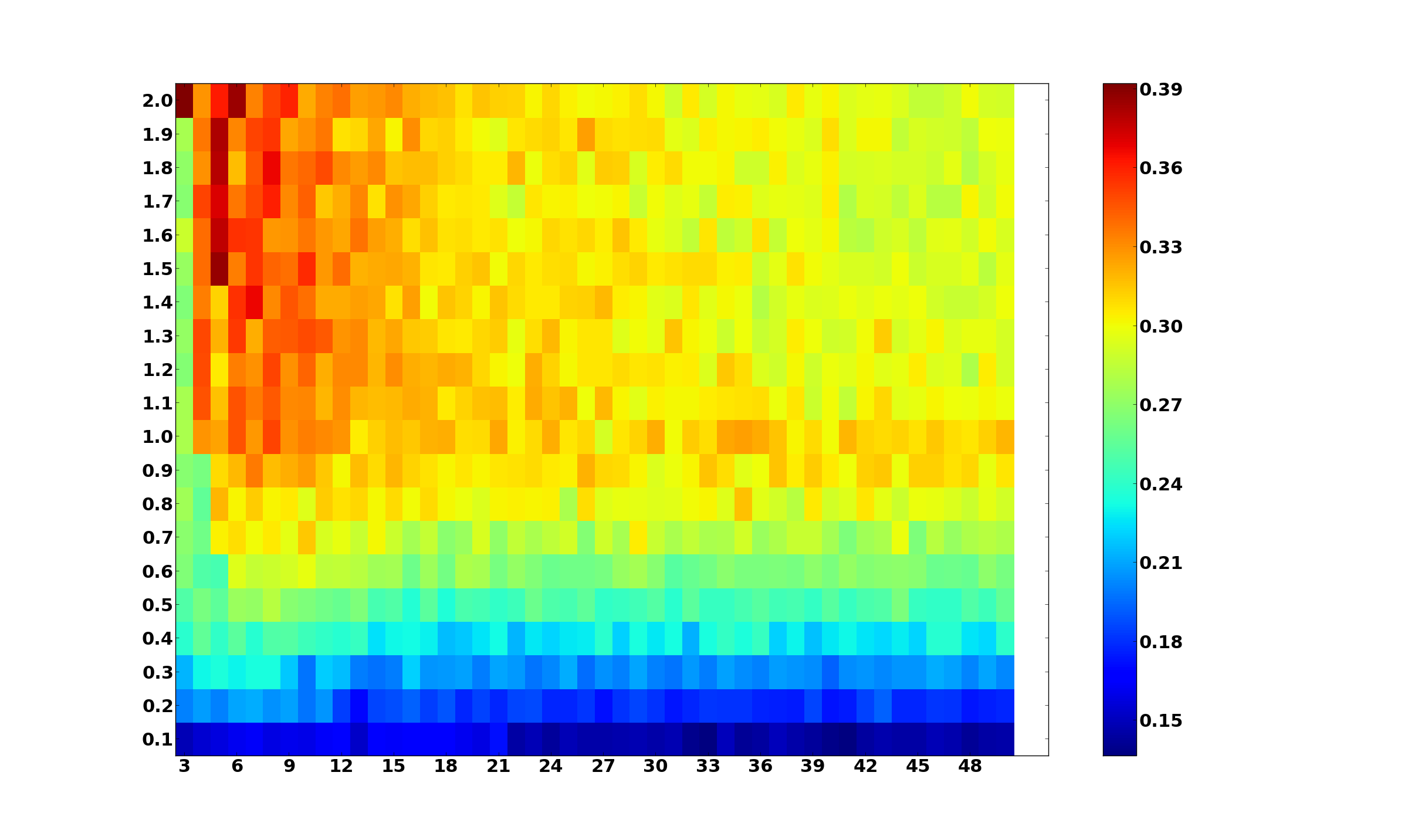}
            \caption{Bayes Mode}
        \end{subfigure}        
        \caption{Separated Moran Process, Sample Size 20. Top Row: Color is given by the inferred values of r minus the true value for Inverse Counting, Counting, and Bayesian Mode (left to right) for 1000 simulations at each pair $(N, r)$. Bottom: Standard deviations of inferred values for same methods. Initial point: $(a, N-a)$ where $a = \text{max}\left(1, \frac{N r}{r+1} \right)$. Although the deviation pattern looks similar, the magnitude is very different from counting to either inference estimate.}
        \label{fps_sims_balanced_sample_20}
\end{figure} 
\end{landscape}

\subsection{Variably-sized Populations and Death-Birth updating}

For well-mixed populations, death-birth updating differs little from birth-death updating other than to effective reduce the population size by one. It also makes fixation slightly more likely when the population is in a state adjacent to an absorbing state. For variable size populations, there are multiple effects in play. Replication events when the population is smaller can carry more information but smaller populations mean that fixation is more likely.

Fixation probabilities can be difficult to derive analytically for variably-sized population processes. The population can now fixate in many ways (the process may have many more than two absorbing states, depending on how variable the population is allowed to be). A discussion of the many ways in which to assign probabilities to birth and death in each round of the process would be lengthy; simulations indicate performance vary similar to the fixed-population size case and will be omitted.

\section{Results: Populations on Graphs}

It will not be possible to cover a comprehensive set of graphs, so we will focus on several interesting cases. Note that to make inferences about $r$, it is not necessary for the graph to be connected or for a particular type to be able to fixate in the population. A sample of a non-absorbing trajectory is enough to make an estimate. The only requirements for directed graphs for each vertex to have an outgoing neighbor (if birth-death is the updating process) or for each vertex to have an incoming vertex (if death-birth is the updating process). For an undirected graph we simply require that each vertex have at least one neighbor. Even these requirements can be relaxed if desired.

\subsection{Cycles and $k$-regular graphs}

Population trajectories on a graph depend on both the number of replicators of each type and the manner in which they are initially distributed. Consider a population on a cycle. One initial case would be for all the replicators of type $A$ to be on a semicircle and all the replicators of type $B$ to be on the other semicircle. In this case, only the replication events at the boundaries of either semicircle will alter the population from its initial state. Similarly, suppose the replicators are initially distributed as $A, B, A, B, \ldots$ around the cycle. Every initial replication event will change the population state in this case. Whether or not this favors one type over the other depends on the true value of $r$. Simulations indicate that inferences from populations on cycles are better than those in well-mixed populations. One reason for this is simply that more replication events may occur on average (depending on the initial distribution) than in the well-mixed case because replicators may have a tendency to replace their own types due to the cycle structure (e.g. in the semicircle initial state). This produces longer trajectories, which can yield better estimates.

The directed cycle has a degenerate special case. For death-birth updating on a directed cycle, the estimate of $r$ will be completely dependent on the prior. This is because each vertex has a single incoming neighbor so no fitness proportionate reproduction occurs during birth events, and so yield no information. In this case, the FPS distribution is improper, with $Pr(r) = 0$ everywhere.

A $k$-regular graph is a graph in which each vertex has exactly $k$ edges. A cycle is a 2-regular graph. For connected $k$-regular graphs, death-birth processes have FPS distribution identical to that of populations of size $k$ ($k-1$ if undirected). For birth death processes, the distribution is that of a population of size $N$, where $N$ is the number of vertices. As noted earlier, smaller populations can have more average deviation from the true value of the parameter $r$ and more variance in estimates, so the processes of birth-death and death-birth can have substantially different behaviors on a $k$-regular graph if $k$ and $N$ are significantly different. This is in contrast to a complete graph, in which the two processes are nearly identical (effective size $N$ vs. $N-1$).

\subsection{Star Topology}

\begin{figure}
    \centering
    \includegraphics[width=0.4\textwidth]{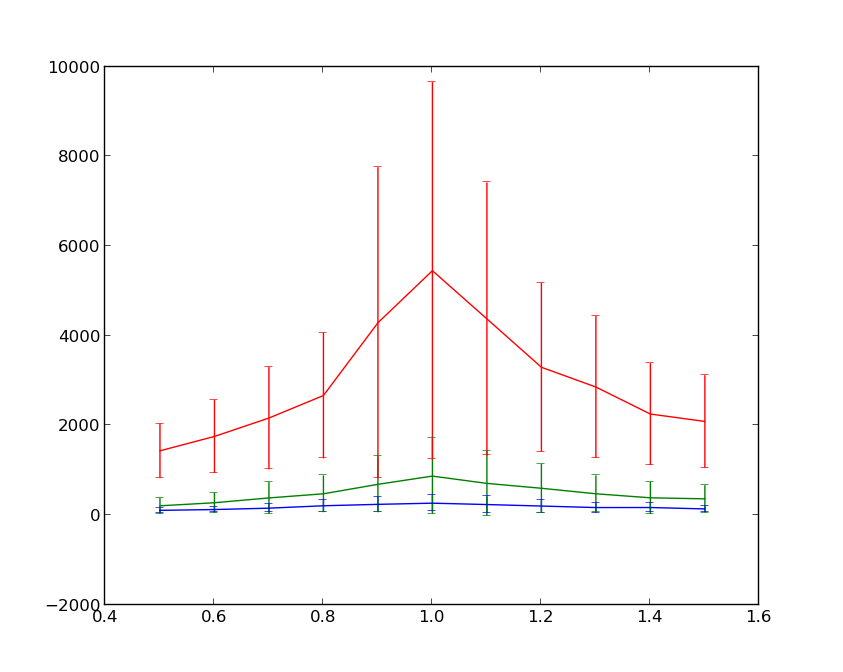}
    \includegraphics[width=0.4\textwidth]{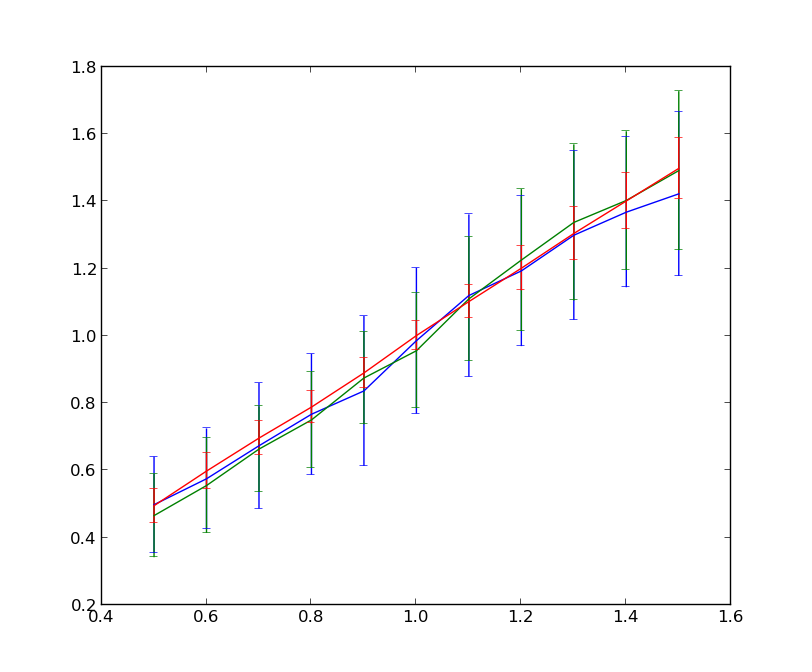}
    \caption{Left: Fitness $r$ vs. Mean fixation time (steps), Population size $N=20$ for complete graph (blue), undirected cycle (green), and star (red). Right: Fitness $r$ vs. Inferred fitness. 200 birth-death simulations per $r$. Initial State: $A$ and $B$ each occupy half of the graph. Longer trajectories produce better estimates with less variance, with the star graph performing best.}
    \label{graph_comp}
\end{figure}

Birth-death processes on star topologies have previously been shown to enhance the strength of selection \cite{lieberman2005evolutionary} versus a well-mixed population. Simulations indicate the star topologies yield more stable inferred values of $r$ in both mean and standard deviation versus well-mixed populations (complete graphs), again likely due to longer trajectories. Again birth-death processes differ significantly from death-birth processes, since in the former case every birth event is from a population of size $N$, and replaces the central vertex with high probability. In the death-birth case, death events occurring on the non-central states are always replaced by the central vertex (which carries no information), and death events at the central vertex act as if the population size were $N-1$. Since death events are equiprobable at every vertex for the death-birth case, the type occupying the central vertex will replicate with probability $(N-1)/N$ regardless of its fitness and increase its proportion in the population! See Figure \ref{graph_comp} for a comparison of the star, cycle, and complete graph.

\subsection{Dynamic Graphs}

It is also possible to infer fitness of birth-death and death-birth processes on graphs with dynamic structure, such as those with active linking \cite{pacheco2006active}. In this case, one would again use the variable-population FPS distribution. As vertices and edges are added or removed, the subpopulations in which fitness proportionate selection occurs may change, and fundamentally little is changed from the case in which an well-mixed population of variable size evolves, from the point of view of computing and estimate from the data. If a vertex has no outgoing neighbors, it replaces itself in the case of a birth-death process. For a death-birth process, there must be at least one incoming neighbor.

\begin{figure}
    \centering
    \includegraphics[scale=0.2]{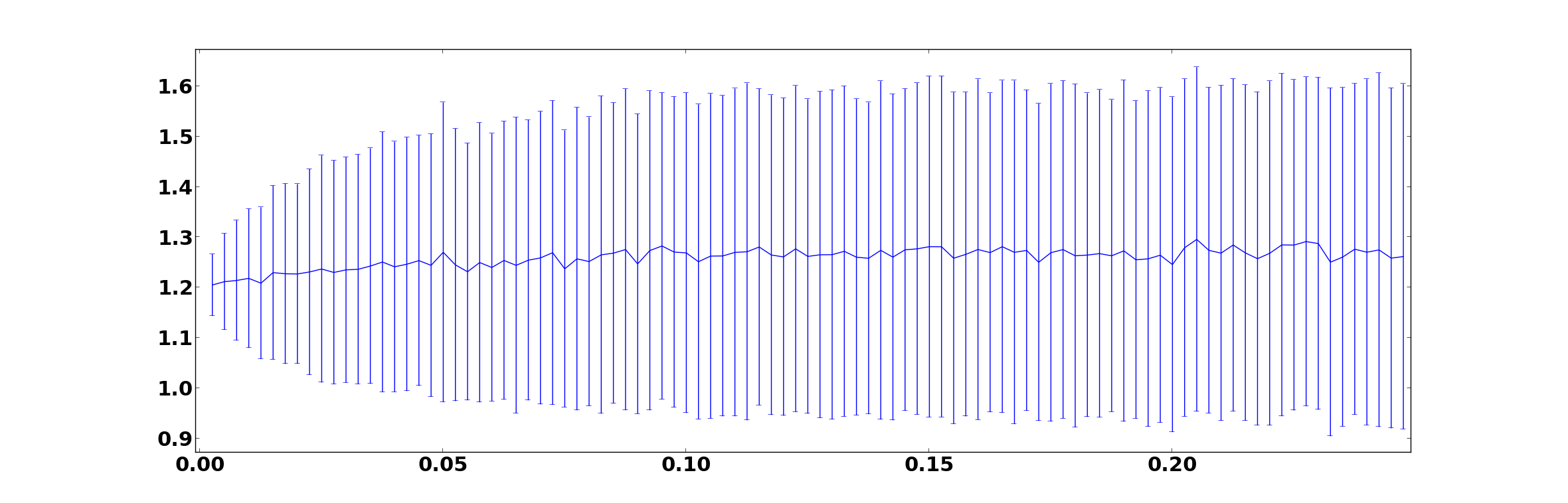}
    \caption{Means and standard deviations of parameter estimates (Bayesian mean) for 1000 birth-death trajectories on a random graph with $p$ on the horizontal axis with $r=1.2$. Initially the population is in the state $(5, 7)$. Outgoing vertices are selected for each iteration (not fixed) with probability $p$. If no outgoing neighbors are randomly selected, the vertex remains occupied in a birth-death event. This can lead to longer trajectories for smaller $p$ and hence better estimates.}
    \label{random_graph}
\end{figure}

The preceding examples of the cycle and the star topology indicate that the effectiveness of inferring the fitness $r$ is affected the population structure, and in particular can improve the estimates obtained. Consider a random graph in which the probability that any two vertices are connected depends on a fixed probability $p$. Means and standard deviations for random graphs for a range of probabilities $p$ are given in Figure \ref{random_graph}. Notice that the estimates are better for small values of $p$, and quickly tend toward a similar distribution as $p$ increases (as the graph becomes ``more complete''). This case is somewhat similar to the case of a $pN$-regular graph. The population gets small-population selection events without the small population tendency to fixate quickly.

\section{Quantifying Information Gain}

\begin{figure}
    \centering
    \includegraphics[scale=0.2]{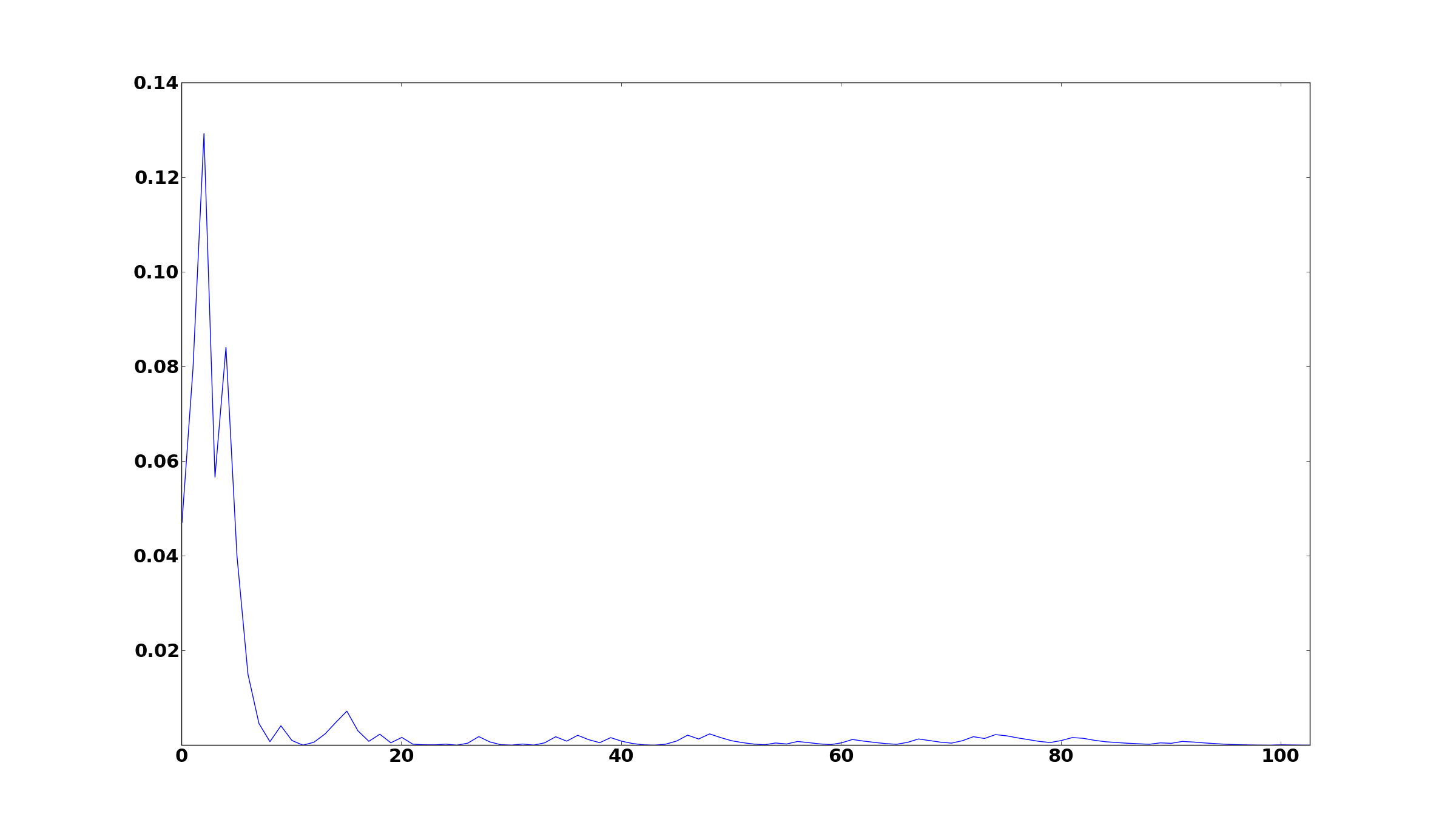}
    \caption{Information gain for a population of size $N=50$, $r=2$, and starting state (33, 17) for a single trajectory.}
    \label{ig_1}
\end{figure}

\begin{figure}
    \centering
    \includegraphics[scale=0.2]{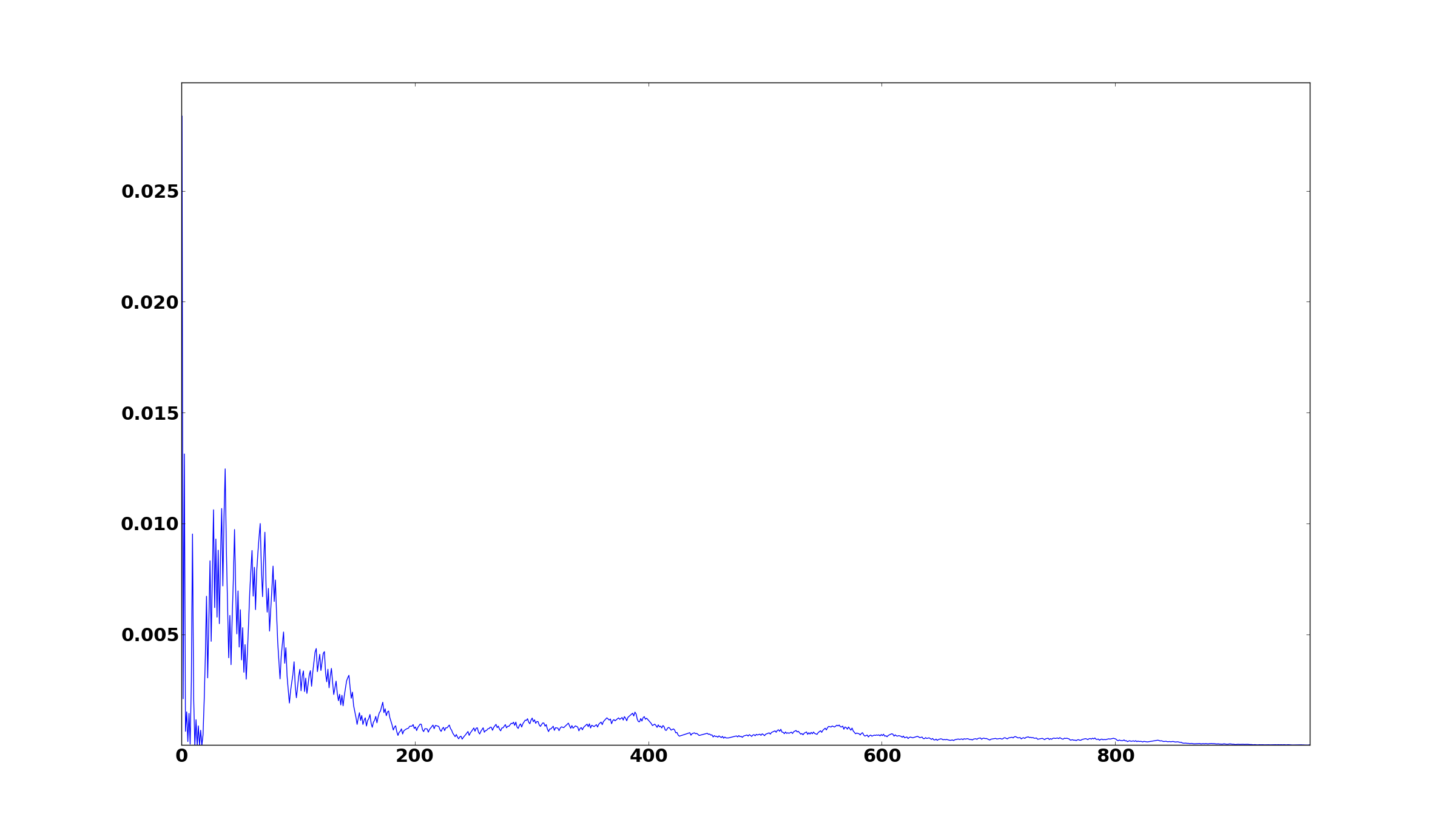}
    \caption{Information gain for a population of size $N=100$, $r=1.2$, and starting state (54, 46) for a single trajectory.}
    \label{ig_2}
\end{figure}

Let us now make good on the promise to quantify how much information is gained by a replication event. Consider the case of a well-mixed fixed size population evolving via the FPS transition probabilities, which in this case form a Bernoulli random variable with regard to the choice of type to replicate. At each step of the trajectory we know the population state $(a, b)$ with $a + b = N$, the estimate of the value of the fitness $r$, and the true value of $r$ (since it was used to create the trajectory). Hence at each step in the process, we can compare the Bernoulli random variables formed by the true value of $r$ and the estimate $\hat{r}$ using the formula $p = a / (a r + b)$ and $\hat{p} = a / (a\hat{r} + b)$. The information divergence between these random variables is 
\[ D_{KL}(p, 1-p; \hat{p}, 1-\hat{p}) = p \log\left(\frac{p}{\hat{p}}\right) + (1 - p) \log \left( \frac{1-p}{1-\hat{p}} \right).\] 
This essentially measures the prediction power of the estimate -- divergences close to zero indicate there is little information left to gain. See Figures \ref{ig_1} and \ref{ig_2} for examples with estimates from Bayesian inference. In both cases, the divergence is near zero after just 20\% of the trajectory data. This partially explains the surprisingly good estimates arising from small samples of trajectories discussed previously.

%
%
%

\section{Discussion}

We have seen that it is possible to fairly accurately infer the unknown fitness of a replicator in Moran-like processes in finite populations on graphs. In particular, even if the trajectory of the population is reduced to a small sample of population state transitions, estimates for the relative fitness can be accurate. Inference is more accurate and has substantially less variation than counting methods, and is much more efficient than using fixation probabilities. Moreover, the method of inference allows the incorporation of prior information, such as a prior belief in the neutral theory of evolution, that is not possible (or not obviously so) with counting methods. In many cases counting methods fail to give a useful estimate because one of the types yields no replication events. This could be worked around by using pseudocounts, adding a nonzero constant $C$ to the numerator and denominator so as to avoid zero/infinite estimates. It is not clear a priori what values would be appropriate choices for $C$, and in any case, a choice of $C$ amounts to an attempt to incorporate prior information into the estimate. Finally, note that Bayesian inference typically yields estimates with far lower variance than the counting methods considered in this work.

At the end of a particular trajectory, the posterior distribution takes on the shape of a normal distribution. This fits the conclusion of the Bernstein von-Mises theorem for Markov processes \cite{borwanker1971bernstein} which says that the posterior should be a normal distribution with mean at the estimate and variance given by the inverse of the Fisher information of the distribution. This has a natural interpretation here. Simply put, from a single trajectory there are many values of the unknown fitness parameter that are likely to produce the trajectory. The posterior distribution indicates the likelihood that the trajectory was generated by any particular value of the fitness (given the prior distribution). Hence the posterior distribution could be used to estimate the probability that the true fitness is greater than one, that it lies in a particular interval, and other similar statistical calculations.

It is possible to infer the fitness of replicators evolving on a structured spacial distribution, such as for evolutionary games on graphs \cite{szabo2007evolutionary}. We saw that the characteristics of the graph alter the effectiveness of inference, with some graphs improving the accuracy and precision of estimates. This may be related to the fact that some graphs can amplify selection \cite{lieberman2005evolutionary}. Dynamic graphs, such as the random graph, behave similarly.

Analogously, it is likely possible to infer fitness of replicating entities in similar dynamical systems. Inferring multiple parameters simultaneously would be of use for Moran processes for more than two types (requiring $n-1$ relative fitness parameters for a $n$-type process), birth-death processes including mutation parameters in addition to (or instead of) fitness parameters, fitness landscapes dependent on game matrices or other additional parameters, Moran-like processes with multiple levels of selection \cite{traulsen2005stochastic}, reproductive processes with mechanisms other than fitness proportionate reproduction. It should also be possible to completely determine fitness landscapes, for instance by observing all the entries of a game matrix on which a fitness landscape is based. All of these examples are straightforward variations of the method described in this manuscript.

\subsection*{Methods}
All computations were performed with python code available at \url{https://github.com/marcharper/mpsim} and \url{https://github.com/marcharper/fitness_inference}. All plots created with \emph{matplotlib} \cite{Hunter:2007}. The software used throughout the results sections depends substantially on the python libraries SciPy \cite{scipy} and Numpy \cite{oliphant2007python}.

\subsection*{Acknowledgments}

This research was supported by the Office of Science (BER),
U. S. Department of Energy, Cooperative Agreement No. DE-FC02-02ER63421.

\bibliography{ref}

\begin{thebibliography}{10}

\bibitem{antal2006fixation}
Tibor Antal and Istvan Scheuring.
\newblock Fixation of strategies for an evolutionary game in finite
  populations.
\newblock {\em Bulletin of mathematical biology}, 68(8):1923--1944, 2006.

\bibitem{borwanker1971bernstein}
J~Borwanker, G~Kallianpur, and BLS~Prakasa Rao.
\newblock The bernstein-von mises theorem for markov processes.
\newblock {\em The Annals of Mathematical Statistics}, pages 1241--1253, 1971.

\bibitem{felsenstein2004inferring}
Joseph Felsenstein and Joseph Felenstein.
\newblock {\em Inferring phylogenies}, volume~2.
\newblock Sinauer Associates Sunderland, 2004.

\bibitem{ficici2000effects}
S~Ficici, J~Pollack, et~al.
\newblock Effects of finite populations on evolutionary stable strategies.
\newblock In {\em Proceedings of the 2000 genetic and evolutionary computation
  conference}, pages 927--934. Morgan-Kaufmann, 2000.

\bibitem{fogel1998instability}
Gary~B Fogel, Peter~C Andrews, and David~B Fogel.
\newblock On the instability of evolutionary stable strategies in small
  populations.
\newblock {\em Ecological Modelling}, 109(3):283--294, 1998.

\bibitem{fudenberg2004stochastic}
Drew Fudenberg, Lorens Imhof, Martin~A Nowak, and Christine Taylor.
\newblock Stochastic evolution as a generalized moran process.
\newblock {\em Unpublished manuscript}, 2004.

\bibitem{hofbauer2003evolutionary}
Josef Hofbauer and Karl Sigmund.
\newblock Evolutionary game dynamics.
\newblock {\em Bulletin of the American Mathematical Society}, 40(4):479, 2003.

\bibitem{Hunter:2007}
J.~D. Hunter.
\newblock Matplotlib: A 2d graphics environment.
\newblock {\em Computing In Science \& Engineering}, 9(3):90--95, 2007.

\bibitem{scipy}
Eric Jones, Travis Oliphant, Pearu Peterson, et~al.
\newblock {SciPy}: Open source scientific tools for {Python}, 2001--.

\bibitem{kimura1985neutral}
Motoo Kimura.
\newblock {\em The neutral theory of molecular evolution}.
\newblock Cambridge University Press, 1985.

\bibitem{kullback1951information}
Solomon Kullback and Richard~A Leibler.
\newblock On information and sufficiency.
\newblock {\em The Annals of Mathematical Statistics}, 22(1):79--86, 1951.

\bibitem{lieberman2005evolutionary}
Erez Lieberman, Christoph Hauert, and Martin~A Nowak.
\newblock Evolutionary dynamics on graphs.
\newblock {\em Nature}, 433(7023):312--316, 2005.

\bibitem{moran1958random}
Patrick Alfred~Pierce Moran.
\newblock Random processes in genetics.
\newblock In {\em Mathematical Proceedings of the Cambridge Philosophical
  Society}, volume~54, pages 60--71. Cambridge Univ Press, 1958.

\bibitem{moran1962statistical}
Patrick Alfred~Pierce Moran et~al.
\newblock The statistical processes of evolutionary theory.
\newblock {\em The statistical processes of evolutionary theory.}, 1962.

\bibitem{nowak2006evolutionary}
Martin~A Nowak.
\newblock {\em Evolutionary dynamics: exploring the equations of life}.
\newblock Belknap Press, 2006.

\bibitem{nowak2004emergence}
Martin~A Nowak, Akira Sasaki, Christine Taylor, and Drew Fudenberg.
\newblock Emergence of cooperation and evolutionary stability in finite
  populations.
\newblock {\em Nature}, 428(6983):646--650, 2004.

\bibitem{oliphant2007python}
Travis~E Oliphant.
\newblock Python for scientific computing.
\newblock {\em Computing in Science \& Engineering}, 9(3):10--20, 2007.

\bibitem{otwinowski2013genotype}
Jakub Otwinowski and Ilya Nemenman.
\newblock Genotype to phenotype mapping and the fitness landscape of the e.
  coli lac promoter.
\newblock {\em PloS one}, 8(5):e61570, 2013.

\bibitem{pacheco2006active}
Jorge~M Pacheco, Arne Traulsen, and Martin~A Nowak.
\newblock Active linking in evolutionary games.
\newblock {\em Journal of theoretical biology}, 243(3):437--443, 2006.

\bibitem{schaffer1988evolutionarily}
Mark~E Schaffer.
\newblock Evolutionarily stable strategies for a finite population and a
  variable contest size.
\newblock {\em Journal of theoretical biology}, 132(4):469--478, 1988.

\bibitem{shaw2010inferring}
Ruth~G Shaw and Charles~J Geyer.
\newblock Inferring fitness landscapes.
\newblock {\em Evolution}, 64(9):2510--2520, 2010.

\bibitem{szabo2007evolutionary}
Gy{\"o}rgy Szab{\'o} and G{\'a}bor F{\'a}th.
\newblock Evolutionary games on graphs.
\newblock {\em Physics Reports}, 446(4):97--216, 2007.

\bibitem{taylor2004evolutionary}
Christine Taylor, Drew Fudenberg, Akira Sasaki, and Martin~A Nowak.
\newblock Evolutionary game dynamics in finite populations.
\newblock {\em Bulletin of mathematical biology}, 66(6):1621--1644, 2004.

\bibitem{traulsen2005stochastic}
Arne Traulsen, Anirvan~M Sengupta, and Martin~A Nowak.
\newblock Stochastic evolutionary dynamics on two levels.
\newblock {\em Journal of theoretical biology}, 235(3):393--401, 2005.

\end{thebibliography}
\bibliographystyle{plain}

\end{document}